\documentclass[12pt,amsfonts,amscd]{amsart}

\usepackage{amssymb,amsmath,amsthm,longtable,multirow,array}
\usepackage{datetime}
\usepackage{amsfonts}
\usepackage{color}
\usepackage{mathrsfs}
\newtheorem{Theorem}{Theorem}[section]
\newtheorem{Proposition}[Theorem]{Proposition}

\newtheorem{Lemma}[Theorem]{Lemma}
\newtheorem{Corollary}[Theorem]{Corollary}
\theoremstyle{definition}
\newtheorem{Definition}[Theorem]{Definition}
\theoremstyle{remark}
\newtheorem{Example}[Theorem]{Example}

\newtheorem*{remark*}{Remark}

\providecommand\ba[1]{\begin{align*}#1\end{align*}}
\providecommand\baa[1]{\begin{align}#1\end{align}}
\providecommand\baaa[1]{\begin{equation}\begin{split}#1\end{split}\end{equation}}

\providecommand\CC{{\mathbb C}}
\providecommand\RR{{\mathbb R}}

\providecommand\NN{{\mathbb N}}
\providecommand\PP{{\mathbb P}}

\providecommand\brs{\begin{remark*}}
\providecommand\ers{\end{remark*}}

\providecommand\eps{\varepsilon}

\providecommand\om{\omega}

\providecommand\Om{\Omega}
\providecommand\al{\alpha}

\providecommand\Gm{\Gamma}

\providecommand\lm{\lambda}

\providecommand\PSH{\operatorname{PSH}}

\providecommand\be{\begin{enumerate}}
\providecommand\ee{\end{enumerate}}
\providecommand\bT{\begin{Theorem}}
\providecommand\eT{\end{Theorem}}
\providecommand\bP{\begin{Proposition}}
\providecommand\eP{\end{Proposition}}
\providecommand\bD{\begin{Definition}}
\providecommand\eD{\end{Definition}}
\providecommand\bE{\begin{Example}}
\providecommand\eE{\end{Example}}
\providecommand\bL{\begin{Lemma}}
\providecommand\eL{\end{Lemma}}
\providecommand\bC{\begin{Corollary}}
\providecommand\eC{\end{Corollary}}

\providecommand\bpp{\begin{proof}}
\providecommand\epp{\end{proof}}
\providecommand\bee{\begin{equation}}
\providecommand\eee{\end{equation}}

\providecommand\Conv{\operatorname{Conv}}

\providecommand\plp{pluripolar\ }
\providecommand\plpt{pluripolar}
\providecommand\psh{plurisubharmonic\ }
\providecommand\psht{plurisubharmonic}

\providecommand\beqq{\begin{eqnarray*}}
\providecommand\eeqq{\end{eqnarray*}}

\DeclareMathOperator{\ddeg}{deg}
\providecommand\vp{\varphi}

\providecommand\iiff{if and only if }

\providecommand\Gm{\Gamma}

\providecommand\ol{\overline}

\providecommand\besp{\begin{split}}
\providecommand\eesp{\end{split}}
\providecommand\bay{\begin{array}}
\providecommand\eay{\end{array}}

\begin{document}
\title{Pluripolar hulls and convergence sets}
\author{Juan Chen and Daowei Ma}

\begin{abstract} The pluripolar hull of a pluripolar set E in $\mathbb{P}^n$ is the intersection of all complete pluripolar sets in $\mathbb{P}^n$ that contain $E$. We prove that the pluripolar hull of each compact pluripolar set in $\mathbb{P}^n$ is $F_\sigma$. The convergence set of a divergent formal power series $f(z_{0}, \dots,z_{n})$ is the set of all ``directions'' $\xi \in\mathbb{P}^{n}$ along which $f$ is convergent. We prove that the union of the pluripolar hulls of a countable collection of compact pluripolar sets in $\mathbb{P}^n$ is the convergence set of some divergent series $f$.  The convergence sets on $\Gamma:=\{[1:z:\psi(z)]: z\in \mathbb{C}\}\subset\mathbb{C}^2\subset\mathbb{P}^2$, where $\psi$ is a transcendental entire holomorphic function, are also studied and we obtain that a subset on $\Gamma$ is a convergence set in $\mathbb{P}^2$ if and only if it is a countable union of compact projectively convex sets, and hence the union of a countable collection of convergence sets on $\Gamma$ is a convergence set.
\end{abstract}

\address{ chen@math.wichita.edu, Department of Mathematics, Wichita State
University, Wichita, KS 67260-0033, USA}
\address{ dma@math.wichita.edu, Department of Mathematics, Wichita State
University, Wichita, KS 67260-0033, USA}
\maketitle

\section{Introduction}

A formal power series $f(z)=f(z_0,z_1,\dots,z_n)=\sum_{\alpha}a_{\al}z^{\al}\in\mathbb{C}[[z_0,z_1,\dots,z_n]]$ with coefficients in $\CC$ is said to be convergent if it converges absolutely in a neighborhood of the origin in $\CC^{n+1}$. Otherwise we say it diverges. So $f$ converges if and only if there is a number $C>0$ such that $|a_{\al}|\le C^{|\alpha|+1}$, for each $\al\in \NN^{n+1}$.  A beautiful classical result of Hartogs (see \cite{Ha}), which can also be interpreted as a formal analog of Hartogs' theorem on separate analyticity, states that a series $f$ converges if and only if $f_{z}(t):=f(z_{0}t,z_{1}t\dots,z_{n}t)$, as a series in $t$, converges for all $z\in \mathbb{C}^{n+1}$. But a divergent power series still may converge in some directions, and this engenders a natural and desirable problem of classifying the convergence sets of divergent power series. Since, for $z\ne 0$, $f_{z}(t)$ converges if and only if $f_{w}(t)$ converges for all $w\in \mathbb{C}^{n+1}$ on the affine line through $z$, we may identify the set of all $z\in \mathbb{C}^{n+1}$ for which $f_{z} $ converges with a subset of $\mathbb{P}^{n}$. The convergence set $\Conv(f)$ of a divergent power series $f$ is defined to be the set of all directions $\xi \in \mathbb{P}^{n}$ such that $f_{z}(t)$ is convergent for some $z\in \pi ^{-1}(\xi )$, where $\pi :\mathbb{C}^{n+1}\backslash \{0\}\rightarrow \mathbb{P}^{n}$ is the natural projection. For the case $n=1$, P.~Lelong (see \cite{Le}) proved that the convergence set of a divergent series $f(z_{1},z_{2})$ is an $F_{\sigma }$ polar set ({\it i.e.}, a countable union of closed sets of vanishing logarithmic capacity) in $\mathbb{P}^{1}$, and moreover, every $F_{\sigma }$ polar subset of $\mathbb{P}^{1}$ is contained in the convergence set of a divergent series $f(z_{1},z_{2})$. The optimal result was later  obtained by A.~Sathaye (see \cite{Sa}) who showed that the class of convergence sets of divergent power series $f(z_{1},z_{2})$ is precisely the class of $F_{\sigma }$ polar sets in $\mathbb{P}^{1}$. In \cite{AM} S.S.~Abhyankar and T.T.~Moh studied convergence sets of divergent power series and showed that a convergence set has Lebesgue measure 0. Levenberg and Molzon (see \cite{LM}) later obtained that a convergence set is \plp and that a compact complete \plp set in $\CC^n$ is a convergence set in $\CC^n$. Convergence sets of formal power series are applied in many areas of mathematics including complex dynamics and investigation of small divisors. In \cite{Pe}  P\'{e}rez-Marco studied polynomial families of dynamical systems presenting problems of small divisors where power series appear naturally and the main result there essentially says that convergence sets in $\CC^n$ are pluripolar. To study the collection $\Conv(\mathbb{P}^{n})$ of convergence sets of divergent power series in higher dimensions, the class $\PSH_{\omega }(\mathbb{P}^{n})$ of $\omega $-plurisubharmonic functions on $\mathbb{P}^{n}$ with respect to the form $\omega:=dd^{c}\log|Z|$ on $\mathbb{P}^{n}$ was considered in \cite{DT}, where the authors proved that $\Conv(\mathbb{P}^{n})$ contains projective hulls of compact pluripolar sets and countable unions of projective varieties, and showed that each convergence set (of divergent power series) is a countable union of projective hulls of compact pluripolar sets. Their main result states that a countable union of closed complete pluripolar sets in $\mathbb{P}^{n}$ belongs to $\Conv(\mathbb{P}^{n})$, which generalized the results of P.~Lelong (see \cite{Le}), Levenberg and Molzon (see \cite{LM}), and Sathaye (see \cite{Sa}).
Convergence sets of formal power series $f(z,t)=\sum_{n=0}^{\infty}f_n(z)t^n$ with $f_n(z)$ holomorphic on a domain $\Omega \subset\CC$ was studied in \cite{ALM}. The authors obtained a characterization of the convergence sets in $\Omega$, which says that a subset of $\CC$ is a convergence set if and only if it is $\sigma$-convex. We get some idea of \cite{ALM} to study convergence sets on $\Gamma$ (see Section 5). 

In the present paper, we consider the class $\PSH_\omega(\PP^n)$ where $\omega$ is the form mentioned above, introduce relative $\omega$-\psh extremal functions and Property J in $\PP^n$, and establish some results on \plp hulls and convergence sets. The main result of this paper generalizes the main result of \cite{DT}.\vspace{0.3cm}\\
\textbf{Theorem~\ref{main}}\, Let $\{K_j\}$ be a sequence of compact pluripolar sets in $\mathbb{P}^{n}$. Then 
\ba{K:=\cup_{j=1}^{\infty}K_{j}^{*}}
 is a convergence set.
\vspace{0.3cm}

Note that $K_j^*$ is the \plp hull of the compact \plp set $K_j$. Of course this is more general than the main result in \cite{DT}, because $K_j^*$ in general is neither compact nor a complete \plp set. In order to prove the main theorem, we introduce the notion of relative $\omega$-\psh extremal functions. We use the relative $\omega$-\psh extremal functions  to decompose the  \plp hull $K_j^*$ into  ``sub-level sets'' $K_j^{(m)}$, and then construct $E_k:=\cup_{j=1}^{k} K_j^{(k)}$ to rewrite $\cup_{j=1}^{\infty} K_j^*$ as a union $\cup_{k=1}^{\infty} E_k$. Unfortunately, in general $E_k$ is not a complete \plp set, and  $E_k$ does not have Property J ({\it i.e.}, $E_k$ does not have Property J with respect to every point in $\PP^n\setminus E_k$), though  $E_k$ does have Property J with respect to each $X\in \PP^n\setminus \cup_{j=1}^{\infty} K_j^*$. Hence with some additional efforts we are able to modify parts of the reasoning in \cite{DT} to finish the proof here.

The following are the precise contents of this paper.

In Section 2, we introduce {\it relative $\omega$-\psh extremal functions} in $\PP^n$ in terms of  $\omega$-\psh functions, gather some useful properties of them and use them to obtain a characterization of the pluripolar hull of a pluripolar set in $\PP^n$ (see Theorem~\ref{1202} and Theorem~\ref{12051}). 

{\it Property J} in $\PP^n$ is defined by means of homogeneous polynomials in Section 3. We establish a useful formula (see Proposition~\ref{11301})
\ba{\sigma(Z,s,K)=\log Q_{K,Z}(e^{-s}),}
which is motivated by \cite[Theorem~5.1.6(i)]{Kl} and \cite[Proposition~2.10]{Sc3}. This formula gives an interpretation of the pluripolar hull  of a \plp set in $\PP^n$ in terms of homogeneous polynomials and hence implies that the pluripolar hull of each compact pluripolar set in $\PP^n$ is $F_\sigma$ (see Theorem~\ref{12052}).  We further establish the connection of Property J, complete pluripolar sets and relative $\omega$-\psh extremal functions in $\PP^n$(see Theorem~\ref{11304}). 

In Section 4, we study convergence sets in $\PP^n$ and prove the main result of the present article that the union of the pluripolar hulls of a countable collection of compact pluripolar sets in $\PP^n$ belongs to $\Conv(\PP^n)$. At the end of Section 4, we consider convergence sets in $\CC^n$ and conclude that the union of the pluripolar hulls of a countable  collection of closed pluripolar sets in $\CC^n$ is a convergence set in $\CC^n$ (see Theorem~\ref{0301}). 

In Section $5$, we consider convergence sets on $\Gamma:=\{[1:z:\psi(z)]: z\in \CC\}\subset\CC^2\subset\PP^2$ with a transcendental entire holomorphic function $\psi$. We prove that a subset on $\Gamma$ is a convergence set in $\PP^n$ if and only if it is a countable union of compact projectively convex sets (see Theorem~\ref{1017}), which implies that the union of a countable collection of convergence sets on $\Gamma$ is a convergence set on $\Gamma$ (see Corollary~\ref{0306}).

\section{Relative $\omega$-\psh extremal functions
} 

It follows from the maximum principle that there are no globally defined non-constant \psh functions on $\PP^n$. However, for a fixed positive closed form $\om$ of bidegree $(1,1)$, the class $\PSH_\omega(\PP^n)$ of $\omega$-\psh functions on $\PP^n$, to be defined below, serves as a substitute of \psh functions. The complete \plp sets on $\PP^n$ can be defined in terms of $\omega$-\psh functions. We introduce relative $\omega$-\psh extremal functions on $\PP^n$, as a substitute of the extremal functions in the local theory (see \cite{BT2}). 

We fix a K\"ahler form 
$$\omega:=dd^c\log|Z|=\frac{i}{2\pi}\partial\overline\partial\log(|Z_0|^2+\cdots +|Z_n|^2)$$
on $\PP^n$, where $d^c=(i/2\pi)(\overline\partial-\partial)$. Note that $\omega$ is the Fubini-Study form on $\PP^{n}$ (see, {\it e.g.}, \cite[p.~30]{GH}).

\bD Let $\Omega\subset \PP^n$ be an open subset of $\PP^n$. A function $\vp:\Omega \to \RR\cup\{-\infty\}$ is said to be {\it $\omega$-\psh} if $\vp+\psi$ is \psh  on $U$ for each open set $U\subset \Om$ and each $C^\infty$ function $\psi$ on $U$ with $dd^c \psi=\om$. 
\eD

For convenience the function 
$\vp\equiv-\infty$ is considered to be both \psh and $\omega$-\psht.

\begin{remark*} (a) Equivalently, an upper semicontinuous function $\vp$ from $\Omega$ to $\RR\cup\{-\infty\}$ is called $\omega$-\psh if $dd^c \vp+\omega\ge0$ (see, {\it e.g.}, \cite{GZ}).  

(b) Observe that $\vp:\Omega \to \RR\cup\{-\infty\}$ is $\om$-\psh if for each open set $U\subset\Omega$ on which there is a $C^{\infty}$ function $\psi$ with $dd^c \psi=\om$ the function $\vp+\psi$ is \psh on $U$. To see this, let us assume that $\psi, \psi_1\in C^\infty(U)$ satisfy $dd^c \psi=\om$ and $dd^c \psi_1=\om$ and that $\vp+\psi$ is \psh on $U$. Then $dd^c (\psi_1-\psi)=0$ and hence $\psi_1-\psi$ is \psh on $U$. It follows that $\vp+\psi_1=(\vp+\psi)+(\psi_1-\psi)$ is \psht.
\end{remark*}

Denote by $\PSH_{\om}(\Omega)$ the family of $\omega$-\psh functions on $\Omega$.
For a homogeneous polynomial $p$, the function $Z\mapsto \log (|p(Z)|^{1/\deg p}/|Z|)$ is a prototypical function in $\PSH_{\om}(\PP^{n})$. Suppose that $\ell(Z):=a_0Z_0+\cdots+a_nZ_n$ is a  linear form and $V$ is an open subset of $U_\ell:=\{Z\in\PP^n:\ell(Z)\ne0\}$. Then a function $\vp$ on $V$ is $\om$-\psh  \iiff the function $\vp(Z)+\log(|Z|/|\ell(Z)|)$ is \psht. Moreover,  a function $\vp$ on $\Om\subset \PP^n$ belongs to $\PSH_\om(\Om)$   \iiff the function $\vp(Z)+\log(|Z|/|\ell(Z)|)$ is \psh on $\Om\cap U_\ell$ for each linear form $\ell$. 

\bD A subset $E$ of $\mathbb{P}^{n}$ is said to be a {\it \plpt} set in $\PP^n$ if there is a function $\vp\in\PSH_\om(\PP^n)$, $\vp\not\equiv -\infty$, such that $E\subset \{\vp=-\infty\}$. A subset $K$ of $\PP^n$ is said to be a {\it complete \plp }set in $\PP^n$ if there is a non-constant function $\vp\in \PSH_\om(\PP^n)$ such that $K=\{\vp=-\infty\}$.\eD

\begin{remark*}
The above definition does not depend on $\omega$; replacing $\omega$ by $c\omega$ with $c>0$ leads to an equivalent definition. In fact, a subset $E$ of $\mathbb{P}^{n}$ is a (complete) \plp set in $\PP^n$ if and only if  the set $U\cap E$ is a (complete) \plp set in $U$ for each affine open set $U$(see, {\it e.g.},  \cite{DT}). 
\end{remark*}

\bD The {\it \plp hull} (see, {\it e.g.},  \cite{LP}) of a pluripolar set E in $\PP^n$ is the intersection of all complete \plp sets in $\PP^n$ that contain $E$, and is denoted by $E^*$, {\it i.e.},
$$E^{*}=\cap\{Z\in \PP^n:\;\; \vp(Z)=-\infty\},$$
where the intersection is taken over all $\om$-\psh functions $\vp$ on $\PP^n$ that are $-\infty$ on $E$.\eD 
\begin{remark*} (a) Notice first that $E\subset E^*$. For the case $n=1$, $E^*=E$. Indeed, since $E$ is polar in $\PP^1$, there exists a $G_\delta$ polar set $F$ such that $E\subset F$. Let $Z\notin E$. By Deny's Theorem (see \break\cite[p.~524]{De}), $F\setminus \{Z\}$ is a complete polar set. 
 It follows that there is a function $\vp\in\PSH_\omega(\PP^1)$ such that $F\setminus\{Z\}=\{\vp=-\infty\}$. Hence, $\vp(Z)>-\infty$, which implies $Z\notin E^*$. Thus $E^*\subset E$, as required.
 
(b) A comment of caution is in order: the \plp hull of a \plp set  is in general not a complete \plp set (see \cite{LP}). However, a result of Zeriahi (see \cite[Proposition~2.1]{Ze}) implies that the \plp hull $E^*$ of a \plp set $E$ is a complete \plp set if $E^*$ is both $G_\delta$ and $F_\sigma$.
\end{remark*}

\bP\label{PlPHp1} If $E_1\subset E_2\subset\PP^n$, then $E_1^*\subset E_2^*$.
\eP 
\bpp This follows immediately from the definition.
\epp 

\bP\label{01012} Let $\{E_j\}$ be a sequence of pluripolar sets in $\mathbb{P}^{n}$. Then
\ba{(\cup_{j=1}^{\infty} E_j)^*=\cup_{j=1}^{\infty} E_j^*.} 
\eP
\bpp  By Proposition~\ref{PlPHp1}, $\cup_{j=1}^{\infty} E_j^*\subset (\cup_{j=1}^{\infty} E_j)^*$. 

Suppose that $Z\notin \cup_{j=1}^{\infty} E_j^*$. Then $Z\notin E_j^*$ for each $j$, and hence there exists $\vp_j\in\PSH_\omega(\PP^n)$ with $\vp_j|_{E_j}=-\infty,\;\vp_j(Z)>-1
$. Set $\vp:=\sum_{j=1}^{\infty}2^{-j}\vp_j$. Then
 \ba{\vp\in \PSH_\omega(\PP^n),\;\;\vp=-\infty \text{\;\;on\;\;} \cup_{j=1}^{\infty} E_j, \text{\;\;and,\;\;} \vp(Z)>-1, } which implies $Z\notin (\cup_{j=1}^{\infty} E_j)^*$. Hence, $(\cup_{j=1}^\infty E_j)^*\subset\cup_{j=1}^\infty E_j^*$. This completes the proof.\epp

\bD\label{01011} For $E\subset\PP^n$, $Z\in\PP^n$, and $s\ge
0$, the {\it relative $\omega$-\psh extremal function}  with respect to the set $E$ and the parameter $s$ is defined by
\ba{\sigma(Z,s,E)=\sup\{\vp(Z): \vp\in\PSH_{\omega}(\PP^n), \vp\le -s\chi_E\},}
where $\chi_E$ is the characteristic function of $E$. Note that $\sigma(Z,0,E)=0$.
\eD 
\begin{remark*} The relative $\omega$-\psh extremal function $\sigma(Z,s,E)$ defined above differs from the relative $\omega$-\psh extremal function in \cite{GZ} and the \psh measure in \cite{BT,Sa2} in that it involves a parameter $s$. Note that $\sigma(Z,s,E)$ depends on $s$ in a nonlinear way, {\it i.e.,} $\sigma(Z,s,E)\ne s\sigma(Z,1,E)$.
In \cite{GZ}, the authors considered $\sigma(Z,s,E)$ for the special case $s=1$.  On the other hand, it is meaningless to involve the parameter $s$ in the corresponding definition in \cite{Sa2}, since there the \psh measure  $\om(z,E)$ was defined in terms of functions belonging to $\PSH(\CC^n)$, which is a convex cone in the sense that if $\alpha$, $\beta$ are non-negative numbers and if $u, v \in \PSH(\CC^n)$ then $\alpha u+\beta v\in\PSH(\CC^n)$, and therefore would satisfy $\om(z,s,E)= s\om(z,1,E)$.
\end{remark*}

For $F\subset\CC^n$, $z\in \CC^n$, and $s>0$, the relative extremal function with respect to the set $F$ and the parameter $s$ is defined by
\ba{\sigma(z,s,F):=\sup\{u(z): u+(1/2)\log(1+|\cdot|^2)\in\PSH(\CC^n), u\le -s\chi_F\}.}

Let $\iota: \CC^n\rightarrow\PP^n$ be defined by 
\baa{\label{IOTA}\iota(z_1,\cdots,z_n)=[1:z_1:\cdots:z_n].}
Then $\iota$ embeds $\CC^n$ into $\PP^n$ and identifies $\CC^n$ with the affine open set 
\ba{\iota(\CC^n)=U_0:=\{Z\in\PP^n: Z_0\ne0\}.}
 Note that $\sigma(\iota(z),s,\iota(F))=\sigma(z,s,F)$.

For $E\subset\PP^n$, let  $\sigma^*(Z,s,E)$ denote the upper semicontinuous regularization of $\sigma(Z,s,E)$ with respect to $Z$. Then $\sigma^*(Z,s,E)$ is $\omega$-\psh in $Z$, which is a straightforward consequence of the analogous local result for sequences of \psh functions (see, {\it e.g.}, \cite{Kl, GZ}).

\bP\label{SMP1} Let $Z\in\PP^n$. If $E_1\subset E_2\subset \PP^n$, then $\sigma(Z,s,E_2)\le \sigma(Z,s,E_1)$ and $\sigma^*(Z,s,E_2)\le \sigma^*(Z,s,E_1)$ for each $Z\in\PP^n$. \eP
Proposition~\ref{SMP1} is a direct consequence of the definition.
\bP\label{SMP2} If $E, P\subset\PP^n$ and if $P$ is pluripolar, then $\sigma^*(Z,s,E)=\sigma^*(Z,s,E\cup P)$. In particular $\sigma^*(Z,s,P)=\sigma^*(Z,s,\varnothing)=0$ when $P$ is \plpt.
\eP
\bpp By Proposition~\ref{SMP1}, it is enough to prove that $\sigma^*(Z,s,E)\le\sigma^*(Z,s,E\cup P)$. Given $\delta>0$, let $\eps>0$ be such that $(s+\delta)/(1+\eps)\ge s$. Let $\vp,\psi\in\PSH_\omega(\PP^n)$ be such that $\vp\le-s\chi_E$, $\psi\le 0$, $\psi\not\equiv -\infty$ and $\psi=-\infty$ on $P$. Then 
\ba{\frac{\vp-\delta+\eps\psi}{1+\eps}\in\PSH_\omega(\PP^n), \;\;\;\; \frac{\vp-\delta+\eps\psi}{1+\eps}\le -s\chi_{E\cup P},}
 which implies that 
\ba{\frac{\vp-\delta+\eps\psi}{1+\eps}\le\sigma^*(\cdot,s,E\cup P).}
Letting $\eps\to 0$, we obtain that $\vp(Z)-\delta\le\sigma^*(Z,s,E\cup P)$ for $Z\in\PP^n\setminus W$ with $W:=\{\psi=-\infty\}$. Note that $W$ has Lebesgue measure zero. It follows that, for each $Z\in\PP^n$,
\ba{\vp(Z)-\delta\le\sigma^*(Z,s,E\cup P),}
which implies that
\ba{\sigma^*(Z,s,E)-\delta\le\sigma^*(Z,s,E\cup P).} 
Letting $\delta\to0$ yields that $\sigma^*(Z,s,E)\le \sigma^*(Z,s,E\cup P)$ , which completes the proof.
\epp

\bD\label{DTH}(see, {\it e.g.}, \cite{Kl}) Let $F$ be a subset of $\CC^n$ and let $\zeta\in\CC^n$. Then $F$ is {\it non-plurithin}  (or {\it non-thin} for the case $n=1$) at $\zeta$ if $\zeta\in \overline{F\setminus \{\zeta\}}$ and if, for each \psh function $u$ defined on a neighborhood of $\zeta$,
\ba{u(\zeta)=\limsup_{z\to\zeta, z\in F\setminus \{\zeta\}}u(z).}
Otherwise we say that $F$ is {\it plurithin} (or {\it thin} for the case $n=1$) at $\zeta$. 
\eD
\begin{remark*} Let $u$ be a \psh function on a neighborhood of $\zeta\in\CC^n$. We certainly have 
$\limsup_{z\to\zeta}u(z)= u(\zeta)$. It follows that a set $F$ is non-plurithin at each point of its interior. In particular, an open set is non-plurithin at each point of itself. 
\end{remark*}

For $E\subset\PP^n$, let \ba{G:=\{Y\in E: \sigma^*(Y,s,E)>-s\}.}
Then the negligible set (see, {\it e.g.}, \cite{Kl}) $G$ is pluripolar. Since $\sigma^*(Z,s,E)\le-s\chi_{E\setminus G}$, we have 
\ba{\sigma^*(Z,s,E)\le\sigma(Z,s,E\setminus G)\le \sigma^*(Z,s,E\setminus G).}  
Thus by Proposition~\ref{SMP2}, 
\baa{\label{SM}\sigma^*(Z,s,E)=\sigma(Z,s,E\setminus G).}

  Denote by $E^{\circ}$ the interior of $E$. 
  
\bP\label{NTP} Let $E\subset\PP^n$ be such that $E^{\circ}$ is non-plurithin at each point of $E$. Then $\sigma^*(Z,s,E)=\sigma(Z,s,E)$ for each $Z\in\PP^n$. \eP
\bpp It suffices to show that $\sigma^*(Z,s,E)\le\sigma(Z,s,E)$ for each $Z\in\PP^n$. Since  $\sigma(\cdot,s,E)=-s$ on $E$, we have $\sigma^*(\cdot,s,E)=-s$ on $E^{\circ}$. By definition if $E^{\circ}$ is non-plurithin at $Z$, then $\sigma^*(Z,s,E)=-s$. Hence,  the assumption that $E^{\circ}$ is non-plurithin at each point of $E$ implies that $\sigma^*(\cdot,s,E)=-s$ on E. It follows that $\sigma^*(Z,s,E)\le\sigma(Z,s,E)$, which completes the proof.
\epp
   
\bC\label{NTP1} If $E\subset\PP^n$ is open, then $\sigma^*(Z,s,E)=\sigma(Z,s,E)$ for each $Z\in\PP^n$.\eC

As an example we now give explicit formulae for $\sigma^*(z, s, \rho B)$, where $z\in\CC^n$, $B=\{z\in\CC^n: |z|\le1\}$, and $\rho B=\{z\in\CC^n: |z|\le\rho\}$. It is an immediate consequence of the definition that $\sigma^{*}(z,0,F)=0$ for each subset $F$ of $\CC^n$. By Proposition~\ref{NTP}, $\sigma^*(z,   s,F)=-s$ for each point $z$ at which the interior $F
^\circ$ is non-plurithin. Thus $\sigma^*(z,s,\rho B)=-s$ for $|z|\le \rho$. So it suffices to calculate $\sigma^*(z,s,\rho B)$ for $s>0$ and $|z|>\rho$.

\bE\label{Sigma}  Let $s>0$ be a fixed positive number and consider\ba{u(z):=\sigma^*(z,s,\rho B)+(1/2)\log(1+|z|^2),\; |z|\ge \rho.} Then the function $u(z)$ is given by  \baa{\label{Sigma1} u(z)=\begin{cases} (1/2)\log(1+\rho^2)-s+\log(|z|/\rho),& \text{ \  if \ }s\ge(1/2)\log(1+\rho^{-2}),\\  (1/2)\log(1+\rho^2)-s+\displaystyle{\frac{\log(|z|/\rho)}{1+e^{-2\lm}}},& \text{ \  if \ }s<(1/2)\log(1+\rho^{-2}),\;\; |z|\le e^\lm,\\ (1/2)\log(1+|z|^2),& \text{ \  if \ }s<(1/2)\log(1+\rho^{-2}),\;\; |z|>e^\lm,\end{cases}} where $\lm=\lm(\rho,s)$ is the unique number in the interval  $(\log\rho,\infty)$ that satisfies the equation \baa{\label{LM}\frac{e^{2\lm}}{1+e^{2\lm}}=(\frac12\log\frac{1+e^{2\lm}}{1+\rho^2}+s)/(\lm-\log\rho).} \eE
\bpp Notice first that $u(z)\in\PSH(\CC^n)$ and for all $z\in \CC^n$,
\ba{u(z)\le \frac12\log(1+|z|^2).} Since 
 \(u(z)= -s+(1/2)\log(1+\rho^{2}) \) on \(\{|z|=\rho\}\), and since  
 the  Siciak's extremal function of $\rho B$ is $\log^+(|z|/\rho)$ (see \cite[p.~185]{Kl}), it follows that  
\ba{u(z)&\le -s+\frac12\log(1+\rho^{-2})+\log|z| \text{\;\;\;\;\;\;for\;\;} |z|\ge\rho.} 
Moreover, since the set $\rho B$ is ``radial'', the function $u$ is radial. 
Set $v(t)=u(e^t,0,\dots,0)$ for \(t\ge \log\rho\). Then $u(z)=v(\log|z|)$. The function $v(t)$ satisfies
\baa{\label{r1}v(t)&\le \frac12\log(1+e^{2t}),\\
       \label{r2}v(t)&\le -s+\frac12\log(1+\rho^{-2})+t. }
Since $u(z)$ is \psht, the function $v(t)$ is increasing and  convex (see \cite[p.~45]{Rt}). Therefore, $v$ is the greatest convex function that does not exceed the function $w(t)$, where 
\ba{w(t):&=\min\{ \mu_1(t), \mu_2(t) \},\\
\mu_1(t):=\frac12\log(1+&e^{2t}),\;\;\mu_2(t):=-s+\frac12\log(1+\rho^{-2})+t.}
Both $\mu_1$ and $\mu_2$ are increasing and convex, so $w(t)$ is increasing; but $w(t)$ is not necessarily  convex.

When $s\ge(1/2)\log(1+\rho^{-2})$, we have $\mu_2(t)\le t<\mu_1(t)$, so $w(t)=\mu_2(t)$ is convex, and therefore $v(t)=w(t)=\mu_2(t)$.
  
We now consider the case $s<(1/2)\log(1+\rho^{-2})$. The function $\mu(t):=\mu_2(t)-\mu_1(t)$ satisfies $\mu(\log\rho)=-s<0$, $\lim_{t\to\infty} \mu(t)=-s+(1/2)\log(1+\rho^{-2})>0$, and $\mu'(t)>0$. It follows that there is a unique point $t_0\in(\log\rho,\infty)$ such that $\mu(t_0)=0$. Thus
\ba{w(t)=\begin{cases} \mu_2(t), & \text{ \  if \ } \log\rho\le t\le t_0, \\ \mu_1(t), & \text{ \  if  \ }  t_0\le t.\\\end{cases} }
Since $w'(t)=1$ on $(\log\rho, t_0)$ and $w'(t)<1$ on $(t_0,\infty)$, it follows that $w$ is not convex on $(\log\rho,\infty)$.
Consider the function 
\ba{\tau(t):=\frac{\mu_1(t)-\mu_2(\log\rho)}{t-\log\rho}, \;\;\; \;\;t\ge t_0,}
which is the slope of the line segment from $(\log\rho, \mu_2(\log\rho))$ to $(t, \mu_1(t))$.  Since $\tau(t_0)=1=\lim_{t\to\infty} \tau(t)$, and since $\tau(t)<1$ for $t>t_0$, we see that  $\tau(t)$ must attain its minimum at some $\lm>t_0$. 
Now $\tau(\lm)$ is necessarily equal to $\mu_1'(\lm)$, and hence the strict convexity of $\mu_1$ implies the uniqueness of $\lm$. Since $\tau(\lm)=\mu_1'(\lm)$, the unique number $\lm$ is determined by (\ref{LM}).  Consequently, the function $v(t)$ is given by
\ba{v(t)=\begin{cases} \mu_2(\log\rho)+\tau(\lm)(t-\log\rho), & \text{ \  if \ } \log\rho\le t\le \lm, \\ \mu_1(t), & \text{ \  if  \ }  \lm\le t.\\\end{cases} }
In summary, combining the two cases, we obtain
 \ba{ v(t)=\begin{cases}
(1/2)\log(1+\rho^{2})-s+t-\log\rho, & \text{if }s\ge(1/2)\log(1+\rho^{-2}), \\ 
(1/2)\log(1+\rho^{2})-s+(1+e^{-2\lm})^{-1}(t-\log\rho), & \text{if }s<(1/2)\log(1+\rho^{-2}),\;\; t\le\lm,\\
(1/2)\log(1+e^{2t}), & \text{if }s<(1/2)\log(1+\rho^{-2}),\;\; t>\lm,\end{cases} }
where $\lm$ satisfies (\ref{LM}). The proof of (\ref{Sigma1}) is complete. This gives explicit formulas for $\sigma^*(z,s,\rho B)$ because it is equal to $u(z)-(1/2)\log(1+|z|^2)$.
\epp
\begin{remark*} By the above implicit expression of the unique number $\lm(\rho,s)$ which depends on $\rho$ and $s$,  we obtain that 
\baa{\label{DLm}\frac{\partial\lm}{\partial\rho}=\frac{1+e^{2\lm}}{2\rho e^{2\lm}(1+\rho^2)}\cdot \frac{e^{2\lm}-\rho^2}{\lm-\log\rho}>0.}
Hence $\lm(\rho,s)$ is increasing with respect to $\rho$, {\it i.e.}, $\lm(\rho_1,s)<\lm(\rho_2,s)$ if $\rho_1<\rho_2$.
\end{remark*}

\bP\label{SMP3} If $\{E_j\}$ is an ascending sequence of subsets in $\PP^n$ and $E=\cup_{j=1}^{\infty}E_j$, then $\sigma^*(Z,s,E_j)$ decreases towards $\sigma^*(Z,s,E)$ for each $Z\in\PP^n$.
\eP
\bpp By Proposition~\ref{SMP2}, $\sigma^*(Z,s,E_j)=0$ if $E_j$ is pluripolar. Thus we may assume that $E_1$ is not \plpt. 
The deceasing sequence $\{\sigma^*(Z,s,E_j)\}$ of $\om$-\psh functions  has an $\om$-\psh limit  $\lim_{j\to\infty} \sigma^*(Z,s,E_j)$. Set $\vp:=\lim_{j\to\infty} \sigma^*(Z,s,E_j)$. Then $\vp\in\PSH_\omega(\PP^n)$ and $\sigma^*(Z,s,E)\le\vp(Z)$. We now prove the reversed inequality $\vp(Z)\le \sigma^*(Z,s,E)$. Put \ba{G:=\cup_{j=1}^{\infty}\{Y\in E_j: \sigma^*(Y,s,E_j)>-s\}.} Then G is \plpt. Note that $\vp= -s$ on $E\setminus G$, which implies that 
\ba{\vp(Z)\le \sigma(Z,s,E\setminus G).}
By (\ref{SM}), 
\ba{\vp(Z)\le \sigma^*(Z,s,E). }
Therefore, $\lim_{j\to\infty} \sigma^*(Z,s,E_j)= \sigma^*(Z,s,E)$. This completes the proof.
\epp

\bP\label{SMP4} If $\{K_j\}$ is a descending sequence of compact subsets in $\PP^n$ and $K=\cap_{j=1}^{\infty}K_j$, then $\sigma(Z,s,K_j)$ increases towards $\sigma(Z,s,K)$ for each $Z\in\PP^n$.
\eP
\bpp Noting that $\lim_{j\to\infty}\sigma(Z,s,K_j)$ exists and that $\lim_{j\to\infty}\sigma(Z,s,K_j)\le\sigma(Z,s,K)$, it suffices to show that $\sigma(Z,s,K)\le\lim_{j\to\infty}\sigma(Z,s,K_j)$. Let $\vp\in\PSH_\omega(\PP^n)$ be such that $\vp\le-s\chi_K$. Given $\eps>0$. By the compactness of $K_j$, $K_j$ is contained in the open set $\{Z\in\PP^n: \vp(Z)<-s+\eps\}$ for sufficiently large $j$. Thus, for such values of $j$, 
\ba{\vp(Z)-\eps\le \sigma(Z,s,K_j)\le \lim_{j\rightarrow \infty} \sigma(Z,s,K_j),  \;\;\;Z\in\PP^n.}
 Consequently, $\sigma(Z,s,K)\le \lim_{j\rightarrow \infty} \sigma(Z,s,K_j)+\eps$. Letting $\eps\to0$ yields $\sigma(Z,s,K)\le\lim_{j\to\infty}\sigma(Z,s,K_j)$, as desired.
\epp

\begin{remark*} (a) For Proposition~\ref{SMP3}, we do not know whether $\sigma(Z,s,E_j)$ decreases towards $\sigma(Z,s,E)$ when $n\ge2$, but we have the following Proposition~\ref{0502} for the case $n=1$. 

(b) For Proposition~\ref{SMP4}, the following Example~\ref{04172} implies that $\lim_{j\to\infty}\sigma^*(Z,s,K_j)\ne\sigma^*(Z,s,K)$ in general. 
\end{remark*}

\bP\label{0502} Let $\{E_j\}$ be an ascending sequence of subsets of $\PP^1$, and let $E:=\cup_{j=1}^{\infty} E_j$. Then \ba{\lim_{j\rightarrow\infty}\sigma(Z,s,E_j)=\sigma(Z,s,E),\;\;\;\;\;Z\in\PP^1.}
\eP 

\bpp If $Z\in E$, then there exists some $j_0$ such that $Z\in E_j$ for each $j>j_0$, and hence $\lim_{j\rightarrow\infty}\sigma(Z,s,E_j)=\sigma(Z,s,E)=-s$. 

Now suppose that $Z\notin E$. It is enough to show that $\lim_{j\rightarrow\infty}\sigma(Z,s,E_j)\le\sigma(Z,s,E)$. Put
\ba{G:=\{Y\in E: \sigma^*(Y,s,E)>-s\}.}
Then $G$ is polar. 
We claim that
\ba{\sigma(Z,s,E)=\sigma(Z,s,E\setminus G).}
In fact, by Proposition~\ref{SMP1},  we have $\sigma(Z,s,E)\le\sigma(Z,s,E\setminus G)$. Take $\eps>0$. Let $\vp\in\PSH_\omega(\PP^1)$ be such that 
$\vp\le-s\chi_{E\setminus G}$, $\vp(Z)>\sigma(Z,s,E\setminus G)-\eps$. Since $G$ is polar, Deny's Theorem (see \cite[p.~524]{De}) provides a function $\psi\in\PSH_\omega(\PP^1)$  such that $\psi\le -s$ on $\PP^1$, $\psi=-\infty$ on $G$ and $\psi(Z)>-\infty$. Set $a:=\sigma(Z,s,E\setminus G)$, $b:=\psi(Z)$, and put $\vp_{\eps}=(1-\eps)\vp+\eps\psi$. Then $\vp_\eps\in\PSH_\omega(\PP^1)$, $\vp_\eps\le0$ on $\PP^1$, and $\vp_\eps\le-s$ on $E$. Hence, $\sigma(Z,s,E)\ge \vp_{\eps}(Z)>(1-\eps)(a-\eps)+\eps b$. Letting $\eps\to0$ yields $\sigma(Z,s,E)\ge a$, which shows $\sigma(Z,s,E)=\sigma(Z,s,E\setminus G)$.\\
Therefore, in view of (\ref{SM}),
\ba{\lim_{j\rightarrow\infty}\sigma(Z,s,E_j)\le\lim_{j\rightarrow\infty}\sigma^*(Z,s,E_j)=\sigma^*(Z,s,E)=\sigma(Z,s,E\setminus G)=\sigma(Z,s,E),}
as required.
\epp

\bE\label{04172} Let $E_k:=(1/k)B$. Then by Proposition~\ref{NTP} for each $k$, $\sigma^{*}(z,s,E_k)=\sigma(z,s,E_k)$ since $E_k^{\circ}$ is non-plurithin at each point of $E_k$. Fixing $s$, by Example ~\ref{Sigma} and (\ref{DLm}), we have $\lm(k,s)<\lm(j,s)$ when $k>j$. Hence we may claim that $\lim_{k\rightarrow\infty} \lm(k,s)=-\infty$. 
In fact, if $\lim_{k\rightarrow\infty} \lm(k,s)=\lm_0\ne-\infty$, then 
\ba{0<\frac{e^{2\lm_0}}{1+e^{2\lm_0}}=\lim_{k\rightarrow\infty} \frac{\frac12\log(1+e^{2\lm_0})+s}{\lm_0-\log(1/k)}=0,}
 a contradiction. It follows that
\begin{equation*}
\lim_{k\rightarrow\infty}\sigma^{*}(z,s,(1/k)B)=\left\{ 
\begin{array}{ll}
0, & \text{ \  for \ } z\neq0, \\ 
-s,  & \text{ \  for  \ }  z=0.\\
\end{array}
\right. 
\end{equation*}

On the other hand, we have that $\sigma(0,s,\{0\})=-s$ and $\sigma(z,s,\{0\})\le 0$. For $z\ne0$, let $\vp(z):=\log|z|-(1/2)\log(1+|z|^2)$. Then $\vp\in\PSH_\omega(\PP^n)$, $\vp\le 0$ on $\PP^n$ and $\vp(0)=-\infty$. Given $0<\delta<1$, we have $\delta\vp\in\PSH_\omega(\PP^n)$, 
\ba{\delta\vp\le0 \text{\;\;on\;\;} \PP^n,\;\;\delta\vp(0)\le-s,}
which implies that $\delta\vp(z)\le \sigma(z,s,\{0\})$. Letting $\delta\rightarrow 0$ yields $0\le\sigma(z,s,\{0\})$ for $z\ne0$. Hence,
\begin{equation*}
\sigma(z,s,\{0\})=\left\{ 
\begin{array}{ll}
0, & \text{ \  for \ } z\neq0, \\ 
-s,  & \text{ \  for  \ }  z=0.\\
\end{array}
\right. 
\end{equation*}
This implies $\sigma^{*}(z,s,\{0\})=0$.\\
Therefore, 
\ba{\sigma^{*}(z,s,\{0\})\ne \lim_{k\rightarrow\infty}\sigma^{*}(z,s,(1/k)B).} 
\eE

As a function of $s$, $\sigma(Z,s,E)$ has the following property. 
\bP For fixed $Z$ and $
E$, $\sigma(Z,s,E)$ is a decreasing, concave function of $s$.
\eP 
\bpp Fixing $Z$ and $E$, it follows straightforwardly from the definition that $\sigma(Z,s,E)$ is decreasing in $s$.

Let $s_1,s_2\ge
0$ and $\lm\in[0,1]$. Given $\eps>0$, there exist $\vp_1,\vp_2\in\PSH_\omega(\PP^n)$ with $\vp_1\le -s_1\chi_E$ and $\vp_2\le -s_2\chi_E$ such that
\ba{\vp_1(Z)>\sigma(Z,s_1,E)-\eps\;\;\;\text{and}\;\;\;\vp_2(Z)>\sigma(Z,s_2,E)-\eps.}
Set $\vp:=(1-\lm)\vp_1+\lm\vp_2$. 
Then $\vp\in\PSH_\omega(\PP^n)$. Since
\ba{\vp&\le-((1-\lm)s_1+\lm s_2)\chi_E,\\
\vp(Z)&>(1-\lm)(\sigma(Z,s_1,E)-\eps)+\lm (\sigma(Z,s_2,E)-\eps),}
it follows that
\ba{\sigma(Z,(1-\lm)s_1+\lm s_2,E)\ge \vp(Z)>(1-\lm)\sigma(Z,s_1,E)+\lm \sigma(Z,s_2,E)-\eps.}
Letting $\eps\to0$ yields that
\ba{\sigma(Z,(1-\lm)s_1+\lm s_2,E)\ge (1-\lm)\sigma(Z,s_1,E)+\lm\sigma(Z,s_2,E).}
Therefore, $\sigma(Z,s,E)$ is concave in $s$.
\epp

\bP\label{0217}  Let $E$ be  a subset of\, $\PP^n$ and let $Z\in \PP^n\setminus E$. Then $\sigma(Z,s,E)=0$ for each $s>0$ if and only if there is a positive number $\tau$ such that $\sigma(Z,s,E)\ge-\tau$ for each $s>0$.
\eP 
\bpp It suffices to show the sufficiency. Let $\tau>0$ be such that $\sigma(Z,s,E)\ge-\tau$ for each $s>0$. Let $\eps>0$ be given. Set $\alpha=(\tau+\eps)/\eps$. Since $\sigma(Z,\alpha s,E)\ge-\tau$, we see that there is a function $\vp\in\PSH_{\omega}(\PP^n)$ such that 
\ba{\vp\le-\alpha s\chi_E\;\;\;\text{and}\;\;\;\vp(Z
)>-\tau-\eps.}
Now $\alpha^{-1}\vp\in\PSH_{\omega}(\PP^n)$, $\alpha^{-1}\vp\le -s\chi_E$, and 
\ba{\alpha^{-1}\vp(Z)>\alpha^{-1}(-\tau-\eps)=-\eps.}
It follows that $\sigma(Z,s,E)>-\eps$.
Letting $\eps\to 0$ yields that $\sigma(Z,s,E)\ge0$, which completes the proof.
\epp

Set 
\ba{E^{(0)}:=\{Z\in\PP^n: \sigma(Z,s,E)=0\;\;\;\text{for each}\;\;\;s>0\}.}

\bT\label{1202} Let $E$ be a pluripolar set in $\PP^n$. Then $E^{(0)}$ is the complement of $E^*$.
\eT
\bpp Suppose that $Z\notin E^*$. Then there exists a function $\vp\in\PSH_\omega(\PP^n)$ such that
\ba{\vp\le0,\;\;\; \vp|_E=-\infty, \;\;\; \text{and}\;\; \vp(Z)>-\infty.}
Let $s, \eps>0$ be given. Choose $\delta\in(0,1)$ so that $\delta\vp(Z)>-\eps$. Note that $\delta\vp\in\PSH_\omega(\PP^n)$ and $\delta\vp\le -s\chi_E$, which implies that $\sigma(Z,s,E)>-\eps$. Letting $\eps\to 0$ yields that $\sigma(Z,s,E)=0$. This means that $Z\in E^{(0)}$.

Conversely, assume that $Z\in E^{(0)}$. Then $\sigma(Z,2^j,E)=0$ for each $j\in\NN$. It follows that there is a sequence of functions $\{\vp_j\}\subset\PSH_\om(\PP^n)$ with $\vp_j\le-2^j\chi_E$ and $\vp_j(Z)>-1$. Setting $\vp:=\sum_{j=1}^{\infty}2^{-j}\vp_j\in\PSH_\omega(\PP^n)$, we obtain 
\ba{\vp|_E&\le\sum_{j=1}^{\infty}(-1)=-\infty,\\
\vp(Z)&>-\sum_{j=1}^{\infty}2^{-j}=-1.}
Thus $Z\notin E^*$, which completes the proof.
\epp

The following Corollary can be derived immediately from Proposition~\ref{01012} and Theorem~\ref{1202}. 

\bC\label{01014} Let $\{E_j\}_{j=1}^{k}$ be a sequence of pluripolar sets in $\PP^n$. Then\ba{(\cup_{j=1}^{k} E_j)^{(0)}=\cap_{j=1}^{k}E_{j}^{(0)}.}\eC

 Let $E$ be a pluripolar set in $\PP^n$. For $\mu,\beta>0$, put
\baa{\label{mubeta}E_{\mu,\beta}:=\{Z\in\PP^n: \sigma(Z,s,E)\le \beta-\frac1\mu s, \;\;\text{for each}\;\;s>0\},}
\ba{ E_{\mu,\infty}:=\cup_{\beta>0} E_{\mu,\beta}.}

\bT\label{12051}  $E^*=\cup_{\mu>0} E_{\mu,\infty}$.
\eT
\bpp Suppose that $Z\in E^*$. By Theorem~\ref{1202} there is an $s_0>0$ such that $\beta:=-\sigma(Z,s_0,E)>0$. Let $\mu=s_0/\beta$. We now show that $\sigma(Z,s,E)\le \beta-s/\mu$ for each $s>0$, which implies that $Z\in E_{\mu,\beta}$. If $s\le s_0$, then $\sigma(Z,s,E)\le0=\beta-s_0/\mu\le \beta-s/\mu$. Now we assume that $s>s_0$. Set $\lm=s_0/s$. By the concavity of $\sigma(Z,s,E)$ we have
\ba{\sigma(Z,s,E)
&=(1/\lm)(\lm \sigma(Z,s,E)+(1-\lambda)\sigma(Z,0,E))\\
&\le (1/\lm)\sigma(Z,\lm s,E)\\
&=-\beta/\lm=-s/\mu.}
Thus \(\sigma(Z,s,E)<\beta-s/\mu\), which means $Z\in E_{\mu,\beta}$. Hence $E^*\subset\cup_{\mu>0} E_{\mu,\infty}$.

Conversely, assume that $Z\in E_{\mu,\beta}$ for some $\mu,\beta$. Then $\sigma(Z,s,E)\le \beta-s/\mu$, which implies that $\sigma(Z,s,E)<0$ for large enough $s$. Hence, $Z\in E^*$ by Theorem~\ref{1202}. Therefore, $\cup_{\mu>0} E_{\mu,\infty}\subset  E^*$. This completes the proof.
\epp  
\begin{remark*} For $m\in\NN$, $\cup_{m=1}^{\infty} E_{m,\infty}=\cup_{m=1}^{\infty} E_{m,m}$. Setting $E^{(m)}:=E_{m,m}$, we have $E\subset E^{(m)}$, $E^{(m)}\subset E^{(m+1)}$, and $E^*=\cup_{m=1}^{\infty} E^{(m)}$.
\end{remark*}

\section{Property J in $\PP^{n}$}

Let $\pi :\mathbb{C}^{n+1}\backslash \{0\}\rightarrow \mathbb{P}^{n}$ denote
the standard projection mapping that maps a point $z=(z_{0},z_{1},...,z_{n})\in \mathbb{C}^{n+1}$ to its corresponding homogeneous coordinates $\pi (z)=[z]=[z_{0}:z_{1}:...:z_{n}]\in \mathbb{P}^{n}$.
Suppose that $z=(z_{0},z_{1},...,z_{n})\in \mathbb{C}^{n+1}\setminus \{0\}$ and $Z=[Z_0:Z_1:\dots:Z_n]\in\PP^n$. Then $Z=\pi(z)$ if and only if
$$[Z_0:Z_1:\dots:Z_n]=[z_{0}:z_{1}:...:z_{n}],$$
or, equivalently,
$$z_jZ_k=z_kZ_j,\;\;\text{ for $j,k=0,\dots,n$.}$$

Denote by $\mathcal{H}(\CC^{n+1})$ the family of all functions $u\in\PSH(\CC^{n+1})$ which are non-negative homogeneous, {\it i.e.}, which satisfies $u(\lm z)=|\lm|u(z)$ for all $\lm\in \CC$ and $z\in\CC^{n+1}$.

Let $\mathscr H(\CC^{n+1})$ denote the set of all homogeneous polynomials (including the zero polynomial of degree $-1$) in $n+1$ variables $z_0, z_1,\dots, z_n$ with complex coefficients. For $k\ge 0$, let $\mathscr{H}_k(\mathbb{C}^{n+1})$ denote the set of homogeneous polynomials $p\in\mathscr H(\CC^{n+1})$ such that $p(\lm z)=\lm^kp(z)$ for all $\lm\in \CC$ and $z\in\CC^{n+1}$. So each $\mathscr{H}_k(\mathbb{C}^{n+1})$ is a $\CC$-vector space and $\mathscr H(\CC^{n+1})=\cup_{k=0}^\infty \mathscr H_k(\CC^{n+1})$.

For $p\in \mathscr H_k(\CC^{n+1})$ with $k\ge 1$, $z\in\CC^{n+1}\setminus \{0\}$ and $Z=[z]\in \PP^n$, set
\ba{\langle p(Z)\rangle :=\frac{|p(Z)|^{1/k}}{|Z|}=\frac{|p(z)|^{1/k}}{|z|}.}
For $p\in\mathscr{H}_0(\mathbb{C}^{n+1})$, let $\langle p(Z)\rangle=|p(Z)|$. Note that $\langle p(Z)\rangle$ is independent of the choice of the representative $z$ and is a well-defined function on $\PP^n$. Furthermore, if $m, k>0$ and $p\in\mathscr H_k(\CC^n)$, then $\langle p^m(Z)\rangle=\langle p(Z)\rangle$.

For a set $E\subset \mathbb{P}^{n}$, put 
\begin{equation*}
\langle p\rangle_{E}=\sup_{Z\in E}\langle p(Z)\rangle,\;\;\;\langle p\rangle=\langle p\rangle_{\PP^n}.
\end{equation*}

Property J defined in \cite{DT} for a set in $\CC^n$ is the very property one uses to prove that the set is a convergence set. We now introduce Property J in $\PP^n$, which will be used to establish a connection between the union of the \plp hulls of a countable collections of compact \plp sets and a convergence set in $\PP^n$ (see Theorem~\ref{main}).
 
\bD A subset $E$ of $\PP^n$ is said to have {\it Property J} (in $\PP^n$) with respect to a point $X\in \PP^n\setminus E$ if there is a sequence $\{p_j\}\subset\mathscr{H}(\mathbb{C}^{n+1})$ and $0<\eta<1$ such that
    \baa {\label{ch2}\langle p_j(X)\rangle>\eta,\,\,\,\,\,\langle p_j\rangle\le 1,\,\,\,\,\,  \lim_{j\to\infty}\langle p_j\rangle_{E}=0.}
The set $E$ is said to have {\it Property J}  if $E$ has {\it Property J} with respect to each $X\in \PP^n\setminus E$.
\eD

 There is a one to one correspondence between $\PSH_\omega(\PP^n)$ and $\mathcal{H}(\CC^{n+1})$. Given $h(z)\in\mathcal{H}(\CC^{n+1})$ and $Z=\pi(z)$, the function
 \ba{\vp(Z):=\log h(z)-\log|z|}
belongs to $\PSH_\omega(\PP^n)$, and $h\mapsto \vp$ is a bijection from $\mathcal{H}(\CC^{n+1})$ to $\PSH_\omega(\PP^n)$.

 By \cite[Theorem~5.1.6(i)]{Kl} or \cite[Proposition~2.10]{Sc3} we have the following Lemma.

\bL\label{11021} If $\vp\in C(\PP^n)\cap\PSH_\omega(\PP^n)$ and $\vp>-\infty$, then
\ba{\vp(Z)=\sup\{\log\langle p(Z)\rangle: p\in \mathscr{H}(\mathbb{C}^{n+1}), \log \langle p\rangle\le \vp \;\;\text{on}\;\; \PP^n\},\;\;\;\;\;\;Z\in\PP^n.}
\eL

\bP\label{11301} Let $K$ be a compact set in $\mathbb{P}^{n}$ and $Z\in \mathbb{P}^{n}$. For $0<t\le1$, define
\ba{Q_{K,Z}(t)=\sup\{\langle p(Z)\rangle: p\in \mathscr{H}(\mathbb{C}^{n+1}), \langle p\rangle\le 1,  \langle p\rangle_{K}\le t\}.}
Then $\sigma(Z,s,K)=\log Q_{K,Z}(e^{-s}) $.\
\eP
\bpp  For each $p\in \mathscr{H}(\mathbb{C}^{n+1})$ with $\langle p\rangle\le 1,  \langle p\rangle_{K}\le e^{-s}$, we see that 
\ba{\log\langle p(Z)\rangle\in \PSH_\omega(\PP^n),\;\; \log\langle p\rangle\le 0,\;\;  \log\langle p\rangle_{K}\le -s.} Thus $\log Q_{K,Z} (e^{-s}) \le\sigma(Z,s,K)$.

In order to obtain the reversed inequality, it suffices to show that $\vp(Z)\le \log Q_{K,Z}( e^{-s}) $ for each $\vp\in \PSH_\omega(\PP^n)$ with $ \vp\le -s\chi_K$. Let $\vp$ be such a function with $\vp(Z)>\sigma(Z,s,K)-\eps$ for arbitrary $\eps>0$. By the approximation theorem (see \cite[Theorem~8.1]{GZ}), there exists a sequence $\{\vp_j\}$ in $C^{\infty}\cap\PSH_\omega(\PP^n)\}$ such that $\vp_j>-\infty$ and $\vp_j$ decreases towards $\vp$. 
In view of the Hartogs Lemma (see, {\it e.g.}, \cite[Theorem~2.9.18 ]{Kl}) applied in $\PP^n$ and in the open set $\Omega=\{Z\in\PP^n: \vp(Z)<-s+\eps\}$, we have $\vp_j\le \eps$ on $\PP^n$ and $\vp_j\le-s+2\eps$ on $K$ for all sufficiently large $j$. And since
\ba{\vp_j(Z)-2\eps\ge \vp(Z)-2\eps>\sigma(Z,s,K)-3\eps,}
it follows from Lemma~\ref{11021} that
\ba{\vp_j(Z)-2\eps=\sup\{\log\langle p(Z)\rangle: p\in \mathscr{H}(\mathbb{C}^{n+1}),\log \langle p\rangle\le\vp_j-2\eps \;\;\text{on}\;\;\PP^n\},\;\;\;\;\;\;Z\in\PP^n.}
Thus there exists a $p\in\mathscr{H}(\mathbb{C}^{n+1})$ with $\langle p\rangle\le e^{-\eps}$ on $\PP^n$ and $\langle p\rangle\le e^{-s}$ on $K$
 such that 
\ba{\log\langle p(Z)\rangle>\vp_j(Z)-2\eps-\eps>\sigma(Z,s,K)-4\eps.}
Therefore
\ba{\log Q_{K,Z}( e^{-s}) >\log\langle p(Z)\rangle> \sigma(Z,s,K)-4\eps.}
Letting $\eps\to 0$ yields that $\log Q_{K,Z}( e^{-s}) \ge\sigma(Z,s,K)$, which completes the proof.\epp
\begin{remark*} Proposition~\ref{11301} establishes a useful formula 
$\sigma(Z,s,K)=\log Q_{K,Z}(e^{-s})$ which is motivated by \cite[Proposition 4.2]{HL}. Consider a compact pluripolar set $K\subset\PP^n$, by Proposition~\ref{11301},
\ba{K_{\mu,\beta}=\{Z\in\PP^n: Q_{K,Z}(t)\le e^\beta t^{1/\mu}, \;\; 0<t\le 1\}}
or, equivalently,
\ba{K_{\mu,\beta}=\{Z\in\PP^n: \langle p(Z)\rangle\le e^{\beta}t^{1/\mu} \text{\;if\;}\,t>0, p\in \mathscr H(\mathbb{C}^{n+1}),\, \langle p\rangle\le 1, \langle p\rangle_{K}\le t\}.}
We now prove that 
\baa{\label{100917}K_{\mu,\beta}=\{Z\in\PP^n: \langle p(Z)\rangle\le e^\beta\langle p\rangle^{1-(1/\mu)}\langle p\rangle_{K}^{1/\mu} \text{\;\;for\;\;} p\in \mathscr H(\mathbb{C}^{n+1})\}.}
For convenience, let us denote by $L$ the set on the right side of  (\ref{100917}). We see immediately that $L\subset K_{\mu,\beta}$. Conversely, assume that $p\in \mathscr H(\mathbb{C}^{n+1})$ is a nonzero homogeneous polynomial and  $Z\in K_{\mu,\beta}$. Setting $\langle p\rangle=M$ and $\langle p\rangle_{K}/M=t$, we obtain that $\langle p/M\rangle= 1$ and $\langle p/M\rangle_{K}= t$, and hence  $\langle p(Z)/M\rangle\le e^{\beta}t^{1/\mu}=e^{\beta}\langle p/M\rangle_K^{1/\mu}$. Thus
\ba{\langle p(Z)\rangle\le e^{\beta}M^{1-(1/\mu)} \langle p\rangle_{K}^{1/\mu}=e^{\beta}\langle p\rangle^{1-(1/\mu)}\langle p\rangle_{K}^{1/\mu},}
and therefore $Z\in L$. We have proved that $K_{\mu,\beta}\subset L$, as desired.
\end{remark*}

\bT\label{12052} For each compact pluripolar set $K$ in $\PP^n$, the pluripolar hull $K^{*}$ is an $F_{\sigma}$ set.
\eT
\bpp By Proposition~\ref{11301}, for each $m\in\NN$, 
\ba{K^{(m)}=\{Z\in\PP^n: Q_{K,Z}(t)\le e^mt^{1/m}
, \;\; 0<t\le 1\},}
that is,
\ba{K^{(m)}=\cap\{Z\in\PP^n: \langle p(Z)\rangle\le e^mt^{1/m}\},}
where the intersection is taken over all $0<t\le 1$ and all homogeneous polynomials $p\in \mathscr H(\mathbb{C}^{n+1})$ with $\langle p\rangle\le 1, \langle p\rangle_{K}\le t$. Since the set $\{Z\in\PP^n: \langle p(Z)\rangle\le e^mt^{1/m}\}$ is closed for each $p$ and each $t$, we see that $K^{(m)}$ is compact for each $m$. Therefore, by Theorem~\ref{12051}, $K^*$ is $F_\sigma$.
\epp

\bP \label{SAU} Let $K$ be  a compact set in $\PP^n$ and let $X\in \PP^n\setminus K$. Then K has Property J with respect to $X$ if and only if $\sigma(X,s,K)=0$ for each $s>0$.\eP
\bpp Fix $s>0$. Suppose that $K$ has {\it Property J} with respect to $X$. Then there is a sequence $\{p_j\}\subset\mathscr{H}(\mathbb{C}^{n+1})$ and $0<\eta<1$ such that
    \ba {\langle p_j(X)\rangle>\eta,\,\,\,\,\,\langle p_j\rangle\le 1,\,\,\,\,\,  \lim_{j\to\infty}\langle p_j\rangle_{K}=0,}
which implies that $Q_{K,X}(e^{-s})>\eta$. By Proposition~\ref{11301} and Proposition~\ref{0217}, $\sigma(X,s,K)=0$ for each $s>0$.  

Conversely, given $0<\eta<1$, assume that $\sigma(X,s,K)=0$ for each $s>0$. Then we have $\sigma(X,j\log 2,K)=0$ for each $j\in\NN$. By Proposition~\ref{11301}, $Q_{K,X}(2^{-j})=1$. It follows that there exists a sequence $\{p_j\}\subset \mathscr H(\CC^{n+1})$ with $\langle p_j\rangle\le 1$ and $\langle p_j\rangle_{K}\le 2^{-j}$ such that \ba{\langle p_j(X)\rangle>Q_{K,X}(2^{-j})-(1-\eta)=\eta.}
Thus
    \ba {\langle p_j(X)\rangle>\eta,\,\,\,\,\,\langle p_j\rangle\le 1,\,\,\,\,\,  \lim_{j\to\infty}\langle p_j\rangle_{K}=0,} 
which implies that $K$ has {\it Property J} with respect to $X$.\epp  
\bC\label{505} Let $K$ be a compact set of $\PP^n$ and let $X\in\PP^n\setminus K$. Then $K$ has {\it Property J} with respect to $X$ if and only if for each $0<\eta<1$ and each $\eps>0$, there is a polynomial $q\in \mathscr H(\mathbb{C}^{n+1})$ such that 
\baa {\label{ch3}\langle q(X)\rangle>\eta,\;\;\langle q\rangle\le 1,\;\; \langle q\rangle_{K}\le \eps.} \eC  
\bpp Let $0<\eta<1$ and $\eps>0$ be given. By Proposition~\ref{11301} and Proposition~\ref{SAU}, $K$ has Property J with respect to $X$ if and only if $Q_{K,X}(\eps)=1$ for each $\eps$ which means that there exists a $q\in \mathscr H(\CC^{n+1})$ with $\langle q\rangle\le 1$ and $\langle q\rangle_{K}\le \eps$ such that 
\ba{\langle q(X)\rangle>Q_{K,X}(\eps)-(1-\eta)=\eta,}
as required.
 \epp 

\bT\label{11304} Let $K$ be a compact pluripolar set  of $\PP^n$. Then the following are equivalent:

(a) $K$ has Property J;

(b) $K$ is a complete pluripolar set;

(c) $\sigma(X,s,K)=0$ for each $s>0$ and each $X\in \PP^n\setminus K$.
\eT
\bpp That $(a)\Leftrightarrow(c)$ follows from  Proposition~\ref{SAU}. 
Now we prove that $(b)\Rightarrow(a)$. In fact, suppose (b) holds. Let $X\in \PP^n\setminus K$. 
Given $0<\eta<1$ and $\eps>0$. Since $K$ is a complete pluripolar set, there exists $\vp\in\PSH_\omega(\PP^n)$ such that  $\vp\le0$ on $\PP^n$, $\{\vp=-\infty\}=K$, and $\vp(X)=\log \eta$. By the approximation theorem (see \cite[Theorem~8.1]{GZ}), there exists a sequence $\{\vp_j\}\subset\PSH_\omega(\PP^n)\cap C^\infty$ with $\vp_j\le0$ on $\PP^n$, $-\infty<\vp_j<-\log j$ on $K$ such that $\vp_j(X)\ge\log \eta+\eps/j$. By Lemma~\ref{11021} there is a $p_j\in \mathscr{H}(\CC^{n+1})$ with $\langle p_j\rangle\le e^{\vp_j} \le1$, and $\langle p_j\rangle_K\le e^{\vp_j}|_K<1/j$ such that 
\ba{\langle p_j(X)\rangle>\exp(\vp_j(X)-\frac{\eps}{j})>\eta.}
Hence, \ba{\langle p_j(X)\rangle>\eta,\;\; \langle p_j\rangle\le1,\;\; \lim_{j\to \infty}\langle p_j\rangle_K\le \lim_{j\to \infty}\frac1j=0,}
which implies that $K$ has Property J with respect to $X\in \PP^n\setminus K$. Since $X\in \PP^n\setminus K$ is arbitrary, we conclude that $K$ has Property J.

Then it suffices to show that $(c)\Rightarrow(b)$. Let $X\in \PP^n\setminus K$ be arbitrary. Suppose that $(c)$ holds. Then $X\in K^{(0)}$ and hence $X$ does not belong to the \plp hull $K^*$ of $K$ by Theorem~\ref{1202}. Thus $K^*=K$.  A theorem of Zeriahi (see \cite[Proposition~2.1]{Ze}) states that if a pluripolar set $F$ is both $F_{\sigma}$ and $G_{\delta}$, and if $F^*=F$, then $F$ is a complete pluripolar set. Note that $K$, as a compact subset of $\PP^n$, is both $F_{\sigma}$ and $G_{\delta}$. It follows that $K$ is a complete pluripolar set, which completes the proof.\epp

\section{Convergence sets}

Denote by $\CC[[z]]:=\CC[[z_0, z_1,\dots, z_n]]$ the ring of formal power series  of $n+1$ variables $z_{0},z_{1},\dots ,$ $z_{n}$ with complex coefficients. Let $f$ be such a formal power series with 
\ba{f(z)=f(z_{0},z_{1},\dots, z_{n})=\sum a_{\alpha}z^{\alpha} \in \mathbb{C}[[z_{0},z_{1},\dots, z_{n}]],}
where $\alpha \in \mathbb{N}^{n+1}$. 
The series $f$  is said to converge if it converges absolutely in a
neighborhood of the origin in $\mathbb{C}^{n+1}$. Otherwise we say it diverges. So  $f$  converges if and only if there is a number $C> 0$ such that $| a_{\alpha}|\le C^{|\alpha|+1}$ for each $\alpha \in \mathbb{N}^{n+1}$.
When $f$ diverges, $f$ may converge in some directions. The {\it convergence set} of $f$, denoted by $\Conv(f)$, is the set of $Z\in \mathbb{P}^{n}$ for which $f_{z}(t):=f(z_{0}t,z_{1}t,\dots ,$ $z_{n}t)$, as a series of one variable $t$, converges (absolutely in some neighborhood of $0$) for some (and hence all) $z=(z_{0},z_{1},\dots ,z_{n})\in \pi^{-1}(Z)$. This is well-defined since, when $z\ne 0$, $f_z(t)$ converges if and only if $f_{cz}(t)$ converges for all $c\in \mathbb{C}\setminus\{0
\}$. Thus the series $f$ converges \iiff $\Conv(f)=\PP^n$.

\bD A subset $E \subset \mathbb{P}^{n}$ is said to be a {\it convergence set} in $\PP^n$ if $E=\Conv(f)$ for some divergent power series $f$. Equivalently, $E$ is a convergence set if $E\ne \PP^n$ and $E=\Conv(f)$ for some $f\in \CC[[z_0,\dots,z_n]]$.\eD

Denote by $\Conv(\PP^{n})$ the collection of all convergence sets in $\PP^n$:
$$\Conv(\PP^n):=\{\Conv(f):\;f\in \CC[[z_{0},z_{1},\dots ,z_{n}]],\;\; f \text{\,\,diverges}\}.$$
Let $f\in\CC[[z]]$ be a divergent series. Since
$$f_z(t):=f(z_{0}t,z_{1}t,\dots ,z_{n}t)=\sum_{j=1}^{\infty}p_j(z)t^j,\;\;p_j\in \mathscr{H}_j(\CC^{n+1}),$$
we see that 
\begin{equation*}
\Conv(f)=\{ Z \in \mathbb{P}^{n}:\sup_{j}\langle p_j(Z)\rangle <\infty\}.
\end{equation*}
In fact, we have the following lemma (see \cite{Sa}).
\bL \label{SBO}\label{bddconv} Suppose that $E\subsetneqq \mathbb{P}^{n}$. Then $E\in \Conv(\PP^n)$ \iiff there exists a countable family $\mathscr F$ of non-constant homogeneous polynomials in $\mathscr H(\CC^{n+1})$ such that
\ba{E=\{Z\in\PP^n: \sup_{p\in \mathscr F}\langle p(Z)\rangle <\infty\}.}
\eL

\bP If $E_j\in \Conv(\PP^n)$ for each $j=1,2,\dots,m$, then $\cap_{j=1}^m E_j\in \Conv(\PP^n)$.
\eP
\bpp Suppose that $E_1,E_2,\dots, E_m\in \Conv(\PP^n)$. Then Lemma ~\ref{SBO} implies that there are countable families $\mathscr F_1, \mathscr F_2,\dots, \mathscr F_m,$ of non-constant homogeneous polynomials such that 
\ba{E_j=\{Z\in\PP^n: \sup_{p\in \mathscr F_j}\langle p(Z)\rangle <\infty\},\;\;\;\;j=1,2,\dots, m.}
It follows that 
\ba{\cap_{j=1}^m E_j=\{Z\in\PP^n: \sup_{p\in \mathscr \cup_{j=1}^m \mathscr F_j}\langle p(Z)\rangle <\infty\}.}
Therefore, $\cap_{j=1}^m E_j\in \Conv(\PP^n)$.
\epp

\bP \label{508} Let $K$ be a compact pluripolar set in $\PP^n$. Then $K^{(0)}=K^{(m)(0)}$ for $m\ge1$.\eP
\bpp Since $K\subset K^{(m)}$,
by Proposition ~\ref{PlPHp1} and Theorem~\ref{12051} we have  
\ba{K^{*}\subset (K^{(m)})^{*}\subset K^{**}=K^*.} It follows that $K^*=(K^{(m)})^{*}$. By Theorem~\ref{1202} we obtain\ba{K^{(0)}=K^{(m)(0)},}as required.\epp

\bT\label{main} Let $\{K_j\}$ be a sequence of compact pluripolar sets in $\mathbb{P}^{n}$. Then $K:=\cup_{j=1}^{\infty}K_{j}^{*}$ is a convergence set. \eT
\bpp  Let $E_k=\cup_{j=1}^{k}K_j^{(k)}$. Then ${E_k}$ is a compact pluripolar set in $\mathbb{P}^{n}$ and  by Theorem~\ref{12051}, 
\baaa{\label{ReW}K&=\cup_{j=1}^{\infty}K_{j}^{*}\\&=\cup_{j=1}^{\infty}\cup_{m=1}^{\infty} K_j^{(m)}\\ &= \cup_{k=1}^{\infty}(\cup_{j=1}^{k}\cup_{m=1}^{k} K_j^{(m)})\\&=\cup_{k=1}^{\infty}(\cup_{j=1}^{k}K_j^{(k)})\\&=\cup_{k=1}^{\infty}E_k.}

Assume that $X\in \PP^{n}\setminus K$ and $k\ge 1$. Then $X\notin K_{j}^{*}$ and hence $X\in K_j^{(0)}$ for each $j$. By Corollary \ref{01014} and Proposition~\ref{508}, $X\in E_k^{(0)}=\cap_{j=1}^{k} K_j^{(0)}$, which implies that $E_k$ has Property J with respect to $X$ by Proposition~\ref{SAU}. For given $0<\eta<1$, by Corollary~\ref{505} there is a homogeneous polynomial $p\in \mathscr H(\CC^{n+1})$ such that 
\ba{\langle p(X)\rangle>\eta,\;\;\langle p\rangle\le 1, \;\;\text{ and }\;\;\langle p\rangle_{E_k}<k^{-1/\beta},}
where $\beta=a/b<1$ is a rational  number with $a, b$ being coprime positive integers so that $(\eta/k)^\beta> 1/2$. \\
Let $d=\ddeg p$ and  $x\in \pi^{-1}(X)$.  Define $q\in \mathscr {H}(\CC^{n+1})$ by  
\ba{q(y)=(ky\cdot \frac{\ol x}{|x|})^d=[k|x|^{-1}(y_0\ol x_0+\cdots+y_n\ol x_n)]^d,}
we have $\ddeg q=\ddeg p=d$, and
\ba{ \langle q\rangle&=\sup_{Z\in \PP^n}\langle q(Z)\rangle\\ &=\sup_{0\ne z\in \CC^{n+1}}\frac{|q(z)|^{1/\deg q}}{|z|}\\
&=\sup_{0\ne z\in \CC^{n+1}}\frac{(k/|x|)|z_0\ol x_0+\cdots+z_n\ol x_n|}{|z|}\\
&=\frac{k}{|x|}\sup_{0\ne z\in \CC^{n+1}}\frac{|z_0\ol x_0+\cdots+z_n\ol x_n|}{|z|}\\&=\frac{k}{|x|}\cdot |x|=k.
}
The penultimate step is by the Cauchy-Schwarz inequality. 
It follows that $\langle q\rangle=k=\langle q(X)\rangle$.
Let $h:=p^aq^{b-a}$. Then $h\in \mathscr H(\CC^{n+1})$ and $\ddeg h=\ddeg(p^aq^{b-a})=bd$. For each $Z\in \PP^n$, 
\ba{\langle h(Z)\rangle&=\frac{|h(Z)|^{1/(bd)}}{|Z|}\\
&=\frac{|p^a(Z)q^{b-a}(Z)|^{1/(bd)}}{|Z|}\\
&=\frac{|p(Z)|^{a/(bd)}}{|Z|^{a/b}}\cdot \frac{|q(Z)|^{(b-a)/(bd)}}{|Z|^{(b-a)/b}}\\
&=(\frac{|p(Z)|^{1/d}}{|Z|})^{a/b}(\frac{|q(Z)|^{1/d}}{|Z|})^{(b-a)/b}\\
&=\langle p(Z)\rangle^{\beta}\langle q(Z)\rangle^{1-\beta}.}
Hence,
\ba{\langle h\rangle&=\sup_{Z\in \PP^n}\langle h(Z)\rangle=\sup_{Z\in \PP^n}\langle p(Z)\rangle^{\beta}\langle q(Z)\rangle^{1-\beta}\le 1^{\beta}\cdot k^{1-\beta}\le k,\\
    \langle h\rangle_{E_k}&=\sup_{Y\in E_k}\langle h(Y)\rangle=\sup_{Y\in E_k}\langle p(Y)\rangle^{\beta}\langle q(Y)\rangle^{1-\beta}\le k^{-1}\cdot k^{1-\beta}=k^{-\beta}<1,\\
    \langle h(X)\rangle&=\langle p(X)\rangle^{\beta}\langle q(X)\rangle^{1-\beta}>\eta^\beta \cdot k^{1-\beta}=(\eta/k)^{\beta}\cdot k>k/2.}
To summarize, there is an $h_{X}\in\mathscr H(\CC^{n+1})$ such that
\ba{\langle h_{X}\rangle&\le k,\;\
 \langle h_{X}\rangle_{E_k}<1,
\;\text{ and }\;
\langle h_{X}(X)\rangle>k/2.}

Let $U_X:=\{Z\in\PP^n: \langle h_X(Z)\rangle>k/2\}$. Then $U_X$ is a neighborhood of $X$.  
The open cover 
$\{U_X: X\in \PP^{n}\setminus K\}$ of $\PP^{n}\setminus K$ contains a countable sub-cover $\{U_{X_j}: j=1,2,\dots\}$. 
Put $h_{kj}=h_{X_j}$. We see that the sequence $\{h_{kj}\}_j$ satisfies 

(i) $\langle h_{kj}\rangle_{E_k}<1$,

(ii) $\langle h_{kj}\rangle \le k$,

(iii) $\cup_{j=1}^\infty\{Z\in\PP^{n}:\langle h_{kj}(Z)\rangle >k/2\}\supset \PP^{n}\setminus K$.\\
 Enumerate the countable sequence $\{h_{kj}\}$ to obtain a single sequence $\{p_\ell\}$ and choose $m_\ell$ so that $\deg p_\ell^{m_\ell}$ is increasing. Set $q_\ell:=p_\ell^{m_l}$ and let $f:=\sum q_\ell$. Then the convergence set of $f$ is given by
\ba{\Conv(f)=\{Z\in \PP^n:\sup_{\ell} \langle q_\ell(Z)\rangle<\infty\}=\{Z\in \PP^n:\sup_{k,j}\langle h_{kj}(Z)\rangle<\infty\}.}

Suppose that $Z\in K$. Then there is an $k_0\ge 2$ such that $Z\in E_{k_0}$. For $k\ge k_0$, $\langle h_{kj}(Z)\rangle \le 1$ and for $k< k_0$, $\langle h_{kj}(Z)\rangle\le k\le k_0-1$. Hence $\sup_{k,j} \langle h_{kj}(Z)\rangle\le k_0-1<\infty$, which implies $K\subset \Conv(f)$.

Suppose that $Z\notin K$. Then by (iii), $\sup_{j} \langle h_{kj}(Z)\rangle\ge k/2$. Hence, $\sup_{k,j} \langle h_{kj}(Z)\rangle=\infty$, which implies that $Z\notin \{Z\in \PP^n:\sup_{k,j}\langle h_{kj}(Z)\rangle<\infty\}$, that is, $Z\notin \Conv(f)$. Thus $\Conv(f)\subset K$.
Therefore,
\ba{K= \Conv(f).}
By Lemma~\ref{SBO},  $K$ is a convergence set in $\PP^n$. 
\epp

\begin{remark*} In the proof of Theroem~\ref{main}, we rewrite the union of the pluripolar hulls of a countable collection of compact pluripolar sets as a union $\cup_{k=1}^{\infty}E_k$ by (\ref{ReW}) and obtain that $E_k$ has Property J with respect to $X\in\PP^n\setminus (\cup_{j=1}^{\infty}K_{j}^{*})$ for each $k$. But in general $E_k$ does not have Property J ({\it i.e.}, $E_k$  does not have Property $J$ with respect to every point in $\PP^n\setminus E_k$) and hence it is not  a complete pluripolar set. Moreover, $K_j^*$ may be noncompact,  and may not be a complete pluripolar set. Therefore, our main result, Theroem~\ref{main}, is more general than Theorem 4.18 in \cite{DT}.
\end{remark*}

We now consider convergence sets in $\CC^n$. Denote by $\mathscr P(\CC^n)$ the family of all polynomials (including the zero polynomial) in $n$ variables $z_0,\cdots, z_n$ with complex coefficients. For $k\ge 0$, let $\mathscr{P}_{k}(\mathbb{C}^{n})$ denote the family of polynomials of degree at most $k$ in $n$ variables $ z_1,\dots, z_n$ with complex coefficients in $\CC$. So each $\mathscr{P}_{k}(\mathbb{C}^{n})=\oplus_{j=0}^{k}\mathscr{H}_j(\CC^n)$ is a $\CC$-vector space and in particular $\mathscr{P}_{0}(\mathbb{C}^{n})=\mathscr{H}_0(\CC^n)=\CC$. For $P\in \mathscr{P}_{k}(\mathbb{C}^{n})$ with $k>0$ and a subset $F$ of $\CC^n$, we set
\ba{\langle P(z)\rangle_k :=\frac{|P(z)|^{1/k}}{\sqrt{1+|z|^2}},\;\;\langle P\rangle_{k,F}=\sup_{z\in F} \langle P(z)\rangle _F, \;\;\langle P\rangle_k=\langle P\rangle_{k,\CC^n}  .}
Note that if $P\in \mathscr{P}_{k}(\mathbb{C}^{n})$ and $m$ is a positive integer, then $\langle P^m(z)\rangle_{km}=\langle P(z)\rangle_k$.

Let $\Lambda_n$ be the set of series $g(t,z)=\sum_{k=0}^{\infty}P_k(z)t^k\in\CC[z_1,z_2,\dots,z_n][[t]]$ with $P_k\in \mathscr P_k(\CC^{n})$. Since we have the convention that the degree of the zero polynomial is $-1$, the polynomials $P_k$ are allowed to be $0$. For $g\in\Lambda_n$, let $\Conv(g)$ be the set of $z\in\CC^n$ for which $g(t,z)$ converges as a series of one variable 
$t$:
\ba{\Conv(g):=\{z\in\CC^n: g(t,z)\in\CC\{t\}\}}
Here we denote by $\CC\{z\}$ the ring of all power series $g(z)\in\CC[[z]]$, where $z=(z_1,\dots,z_n)$, that are absolutely convergent in a neighborhood of the origin in $\CC^n$. By Hartog's theorem (see \cite{Ha}), $\Conv(g)=\CC^n$ if and only if $g\in\CC\{t,z_1,\dots,z_n\}$, where $\CC\{t,z_1,\dots,z_n\}$ denote the set of convergent series in $n+1$ variables $t,z_1,\dots,z_n$.
\bD A subset $F$ of $\CC^n$ is said to be a {\it convergence set} in $\CC^n$ if $F=\Conv(g)$ for some divergent power series $g\in\Lambda_n$.\eD 

The following proposition can be found in \cite[Proposition 5.3]{DT}. 
\bP\label{180201} Suppose that $F\subsetneqq \mathbb{C}^{n}$. Then $F$ is a convergence set in $\CC^n$ \iiff $\iota(F)=U_0\cap E$ for some $E\in\Conv(\PP^n)$.\eP

\begin{remark*} 
It was shown in \cite{DT} that for a convergence set $F$ in $\CC^n$, the set $\iota(F)$ may or may not be a convergence set in $\PP^n$.
\end{remark*}

Let $W$ be a compact pluripolar subset of $\CC^n$. Then the \plp hull   of $W$ is given by 
\baa{\label{plphc}W^{*}=\iota^{-1}(U_0\cap(\iota(W))^*).}
It follows that $W^*=\cup_{m=1}^{\infty}W^{(m)}$, where
\ba{W^{(m)}=\{z\in\CC^n: \langle P(z)\rangle\le e^m\langle P\rangle_{k
}^{1-1/m}\langle P\rangle_{k,W}^{1/m} \text{\;\;\;for \;\;} k\in\NN, P\in\mathscr P_k(\CC^n)\}.}

By Theorem ~\ref{main}, Proposition~\ref{180201} and in view of (\ref{plphc}), we have  
\bT\label{0301} The union of the pluripolar hulls of a countable  collection of closed pluripolar sets in $\CC^n$ is a convergence set in $\CC^n$.
\eT

\section{Convergence sets on $\Gamma$}

We fix a non-polynomial entire holomorphic function $\psi(z)$ on the complex plane, say $\psi(z)=e^z$. Denote by $\Gm$ the image under the embedding $\CC^2\to\PP^2$ of the graph of $\psi$: 
\ba{\Gm:=\{[1:z:\psi(z)]\in\PP^2: z\in \CC\}.}
For a subset $S$ of $\CC$, let 
$\Gm(S):=\{[1:z:\psi(z)]\in\PP^2: z\in S\}$. Obviously, every subset of $\Gm$ is $\Gm(S)$ for some $S\subset \CC$.

In this section we obtain a characterization of subsets of $\Gm$ which are convergence sets in $\PP^2$. Such a set is called a {\it convergence set on $\Gm$}. Observe that a convergence set on $\Gm$ is necessarily a convergence set in $\CC^2$. The main idea of this section is motivated by  \cite[Theorem~9.2]{HL} and some tools we use here to obtain the main result of this section come from \cite{ALM}.

\bD(see, {\it e.g.}, \cite{HL}) The {\it projective hull} $\hat K$ of a compact set $K\subset\PP^n$ is the set of all points $Z\in\PP^n$ for which there exists a constant $C=C_Z>0$ such that
\ba{\langle p(Z)\rangle\le C\langle p\rangle_K}
for all homogeneous polynomials $p\in \mathscr H(\CC^{n+1})$. A compact set $K\subset\PP^n$ is said to be {\it projectively convex} if $\hat K=K$.
\eD

Note that each algebraic variety in $\PP^n$ is projectively convex. In particular, each finite set is projectively convex  (see, {\it e.g.}, \cite{DT}).

\begin{remark*} (a) The projective hull of a compact set $K\subset\PP^n$ defined above is equivalent to that in \cite[p.~607]{HL} because for $k\ge 1$, the set $H^0(\PP^n,\mathcal{O}(k))$ of global holomorphic sections of the line bundle $\mathcal{O}(k)$ is canonically identified with the set $\mathscr H_k(\CC^{n+1})$ of homogeneous polynomials of degree $k$. We caution that the above definition for projective hulls differs substantially from that of $\hat X_h$ (the $h$ in the subscript stands for ``homogeneous'') in\break \cite[p.~116]{Sn} in that the above definition contains a constant $C$ while the definition in \cite{Sn} requires that $C=1$. For instance, for 
$X=\{(z,0): z\in\CC, |z|\le 1\}\subset \CC^2\subset \PP^2$, the projective hull  $\hat X$ (in $\CC^2$) is $\{(z,0): z\in\CC\}$, while the $\hat X_h$ in \cite{Sn} equals $X$.

(b) For a compact \plp set $K$ in $\PP^n$, by (\ref{100917}) we have
\ba{K_{1,\beta}=\{Z\in\PP^n: &\langle p(Z)\rangle\le e^{\beta}\langle p\rangle_{K} \text{\;\;for\;\;} p\in \mathscr H(\mathbb{C}^{n+1})\}.}
It follows that $\hat K=K_{1,\infty}=\cup_{\beta>0} K_{1,\beta}$.  Furthermore, by Theorem~\ref{12051} we obtain that $\hat K\subset K^*$, which was also obtained  in \cite[Proposition~3.6]{DT} with a different proof.
\end{remark*}

By \cite[Lemma~4.7]{DT} we have the following Proposition.
\bP\label{09221} Let $\{K_j\}$ be an ascending sequence of projectively convex, compact, pluripolar sets in $\PP^n$ and let $K:=\cup_{j=1}^{\infty} K_j$. Let $\{U_j\}$ be a sequence of open sets such that $K_j\subset U_j$ and $\cup_{m=1}^{\infty}\cap_{j=m}^{\infty}U_j\subset K$. Then $K$ is a convergence set in $\PP^n$.
\eP

\bC\label{9222} Let $\{K_j\}$ be an ascending sequence of projectively convex, compact, pluripolar sets in $\PP^n$ such that $K:=\cup_{j=1}^{\infty} K_j$ is $G_\delta$. Then $K\in \Conv(\PP^n)$.
\eC 

\bD\label{1020} The {\it polynomial hull} $\tilde W$ of a compact set $W\subset \CC^n$ is the set of all points $z\in\CC^n$ such that 
\ba{| p(z)|\le \parallel p\parallel_{W}}
for all polynomials $p\in \mathscr{P}(\mathbb{C}^{n})$. A compact set $W\in\CC^n$ is said to be {\it polynomially convex} if $\tilde W=W$.
\eD

\bP\label{1022} Let $W_1, W_2,\cdots, W_k$ be pairwise disjoint compact sets in $\CC
$. Then $\cup_{j=1}^{k} \tilde{W_j}$ is the polynomial hull of $\cup_{j=1}^{k} W_j$. 
\eP
\bpp It suffices to show that the polynomial hull of $\cup_{j=1}^{k} W_j$ is contained in $\cup_{j=1}^{k} \tilde{W_j}$ for the case $k=2$. Suppose that $W_1$, $W_2$ are compact sets with $W_1\cap W_2=\varnothing$. Let $U$ be a bounded connected component of $\CC\setminus (W_1\cup W_2)$, and let $D=:\partial\tilde{\overline{U}}$ be the boundary of the polynomial hull of the closure of $U$. Then 
\ba{D\subset \partial\overline{U}\subset \partial U\subset \partial (\CC\setminus(W_1\cup W_2))=\partial (W_1\cup W_2)\subset (W_1\cup W_2)^-=W_1\cup W_2. }
Since $D$ is connected, it follows that $D\subset W_1$ or $W_2$, and hence $U\subset \tilde D\subset \tilde{W}_1\cup\tilde{W}_2$. Therefore, the polynomial hull of $W_1\cup W_2$ is contained in $\tilde{W}_1\cup\tilde{W}_2$.
\epp

By \cite[Theorem 9.2]{HL}, which depends on a deep theorem in \cite{Sa1}, we have the following Lemma.
\bL \label{1012}Let $W$ be a compact subset in $\CC$, the set $\Gamma(W)$ is projectively convex in $\PP^2$ if and only if\, $W$ is polynomially convex in $\CC$.
\eL

For $r>0$ and $S\subset \CC$, we define the $r$-neighborhood $N_r(S)$ of $S$ by
\ba{N_r(S)=\{z\in\CC:|z-s|<r \,\,\,\text{for some}\,\,\, s\in S \},}
and denote $d(z,S)$ the Euclidean distance between the point $z$ and the set $S$.

\bP\label{1016} Let $\{W_j\}$ be a sequence of compact, polynomially convex sets in $\CC$ and let $F=\cup_{j=1}^{\infty}W_j$. Then there exists an ascending sequence $\{F_k\}$ of compact, polynomially convex sets in $\CC$ such that 
\baa{\label{6.81}F=\cup_{k=1}^{\infty}F_k,}
and such that the neighborhoods $V_k:=N_{1/(3k)}(F_k)$ satisfy 
\baa{\label{6.82}\cap_{k=m}^{\infty}V_k\subset F}
for each integer $m\ge 1$.
\eP
\bpp For $1\le j\le k$, set
     \ba{L_{kj}=W_j\setminus \cup_{\ell=1}^{j-1} N_{1/k}(W_{\ell}), \text{\;\;and\;\;}K_{kj}=\tilde{L}_{kj}.}
Then we have 
\ba{L_{k\ell}\cap L_{kj}=\varnothing\;\;\;\; \text{for}\;\; 1\le\ell<j\le k,}
and 
\ba{L_{kj}\subset L_{k+1,j}.} 
Let $G_k=\cup_{j=1}^{k}L_{kj}$, and $F_k=\cup_{j=1}^{k}K_{kj}$. Then $G_k\subset G_{k+1}$ and by Proposition~\ref{1022},  $F_k=\tilde{G_k}$ for each $k$. 
Hence $F_k\subset F_{k+1}$ and each $F_k$, as the polynomial hull of $G_k$, is polynomially convex. 

Suppose that $z\in F$. Let $j$ be the least integer such that  $z\in W_j$. If $j=1$, then $z\in W_1=F_1$. Assume $j>1$ and choose $k\ge j$ so that $d(z,\cup_{i=1}^{j-1} W_i)>1/k$. Then 
\ba{z\in L_{kj}\subset K_{kj}\subset F_k,}
which implies $F\subset\cup_{k=1}^{\infty}F_k$.

 Conversely, since 
\ba{F_k=\cup_{j=1}^{k}K_{kj}\subset \cup_{j=1}^{k} W_j\subset F,}
we see that $\cup_{k=1}^{\infty}F_k\subset F$. This completes the proof of (\ref{6.81}). 

 We now prove (\ref{6.82}). Fixing an integer $m\ge 1$, we will show that the assumption $(\cap_{k=m}^{\infty}V_k)\setminus F\neq \varnothing$ leads to a contradiction. Assume that $z\in (\cap_{k=m}^{\infty}V_k)\setminus F$. We now prove by induction that 
\baa{\label{xc2}z\in N_{1/(3k)}(\cup_{j=1}^{m}K_{kj}),\;\; \text{for}\;\; k\ge m.}
When $k=m$, 
\ba{z\in V_m=N_{1/(3m)}(\cup_{j=1}^{m}K_{mj}).} 
We assume that $M>m$ and (\ref{xc2}) holds for $k=M-1$, that is, 
\baa{\label{xc3}z\in N_{1/(3M-3)}(\cup_{j=1}^{m}K_{M-1,j}).}

Set 
\ba{Q=\cup_{j=m+1}^{M}K_{Mj},\;\;\;\;\;
    R=Q\cap(\cup_{j=1}^{m}K_{Mj}),\;\;\;\;\;
    S=(\cup_{j=1}^{m}K_{Mj})\setminus R.}
Then $Q\cap S=\varnothing$, and $Q\cup S=F_M$.

Consider a point $w\in R$. There exist $i,\ell$ with $1\le i\le m<\ell\le M$ such that $w\in K_{Mi}\cap K_{M\ell}$. By the construction of $L_{M\ell}$ and $K_{M\ell}$, $w\in K_{M\ell}\setminus L_{M\ell}$ and hence $w$ belongs to a bounded connected component $U$ of the complement $L^c_{M\ell}$ of $L_{M\ell}$. The connected open set $N_{1/M}(w)$ is contained in $L^c_{M\ell}$, so it is contained in some connected component of $L^c_{M\ell}$; since $w\in U$, we obtain that $N_{1/M}(w)\subset U\subset K_{M\ell}$. Thus, 
\baa{\label{09023} N_{1/M}(R)\subset Q.}

Now we consider a point $w\in S$ and a $j$ with $m+1\le j\le M$. Since $w\in N_{1/M}(\cup_{i=1}^{m} W_i)$ and since $w\notin K_{Mj}$, we see that $w$ belongs to the unbounded connected component $V$ of $L^c_{Mj}$. The connected open set $N_{1/M}(w)$ is contained in $L^c_{Mj}$, so it must be contained in $V$. It follows that $N_{1/M}(w)\cap K_{Mj}=\varnothing$. Therefore, 
\baa{\label{09024} N_{1/M}(S)\cap Q=\varnothing.}
By the induction hypothesis, we have
\ba{z\in N_{1/(3M-3)}(\cup_{j=1}^{m}K_{M-1,j})\subset N_{1/(3M-3)}(\cup_{j=1}^{m}K_{Mj})=N_{1/(3M-3)}(S\cup R)),}
and hence 
\baa{\label{09021} z\in N_{1/(3M-3)}(S)\cup N_{1/(3M-3)}(R).}
However, since $N_{1/(3M-3)}(R)\subset N_{1/M}(R)\subset Q\subset F$ by (\ref{09024}), we see that
\baa{\label{09022} z\notin N_{1/(3M-3)}(R).}
Now (\ref{09021}) and (\ref{09022}) imply that 
\baa{\label{09025}z\in N_{1/(3M-3)}(S).}
Making use of (\ref{09024}) and (\ref{09025}), we obtain that
\ba{d(z,Q)&\ge d(S,Q)-d(z,S)\\ 
          &\ge \frac{1}{M}-\frac{1}{3M-3}\\ &\ge \frac{1}{M}-\frac{2}{3M}\\&=\frac{1}{3M},}
and hence $z\notin N_{1/(3M)}(Q)$. This, together with $z\in V_M=N_{1/(3M)}(S\cup Q)$, implies that  
\ba{z\in N_{1/(3M)}(S)\subset N_{1/(3M)}(\cup_{j=1}^{m}K_{Mj}).} 
This completes the proof of (\ref{xc2}).

Since $\cup_{j=1}^{m}K_{kj}\subset \cup_{j=1}^{m} W_j$, (\ref{xc2}) implies that
\ba{d(z,\cup_{j=1}^{m}W_{j})<\frac{1}{3k}, \;\;\text{for each}\;\; k\ge m. }
Letting $k\rightarrow\infty$ yields that $z\in \cup_{j=1}^{m}W_{j}\subset F$, which is a contradiction, as desired.
\epp

\bT\label{1017} $E\subset \Gamma$ is  a convergence set if and only if $E=\cup_{j=1}^{\infty}E_j$, where $E_j$ is  compact and projectively convex for each $j$.
\eT
\bpp Suppose that $E\subset \Gamma$ is  a convergence set and $E=\Conv(f)$ with $f=\sum_{m=1}^{\infty}h_m$, where $h_m\in\mathscr{H}_{m}(\mathbb{C}^{3})$. Then $E=\cup_{j=1}^{\infty}E_j$, where
\ba{E_j:=\Gamma(W_j), \text{\;\;with\;\;\;}W_j:=\{z\in\CC: 
 \langle h_m(1,z,\psi(z))\rangle\le j, \forall\; m\}.}
Since $W_j$ is compact and polynomially convex in $\CC$ for each $j$,  by Lemma ~\ref{1012}, it follows that each $E_j$ is compact and projectively convex in $\PP^2$. 

Conversely, suppose that $E=\cup_{j=1}^{\infty}E_j$, where $E_j$ are  compact and projectively convex. For a fixed positive integer $j$, let $W_j$ be the subset of $\CC$ with $E_j=\Gamma(W_j)$. By Lemma~\ref{1012}, $W_j$ is compact and polynomially convex in $\CC$. Let $F=\cup_{j=1}^{\infty}W_j$. By Proposition~\ref{1016}, there exists an ascending sequence $\{F_j\}$ of compact and polynomially convex sets such that $F=\cup_{j=1}^{\infty}F_j$, and 
\ba{\cup_{m=1}^{\infty}\cap_{j=m}^{\infty} V_j\subset F,}
where $V_j:=N_{1/(3j)}(F_j)$. Consequently, there exists a sequence $\{U_j\}$ of open sets with $E_j\subset U_j$ on $\Gamma$ such that 
 \ba{\cup_{m=1}^{\infty}\cap_{j=m}^{\infty} U_j\subset E.}
By Proposition~\ref{09221}, $E\in \Conv(\PP^2)$.
\epp

\begin{remark*}
  We know that the infinite intersection of a countable collection of convergence sets is not necessarily a convergence set (see \cite[Proposition~4.15]{DT}). We do not know whether the union of a countable collection of convergence sets is necessarily a convergence set. However, by Theorem ~\ref{1017}, we have the following
\end{remark*}

\bC\label{0306} The union of a countable collection of convergence sets on $\Gamma$ is a convergence set on $\Gamma$.
\eC


\begin{thebibliography}{}

\bibitem{AM} S.S.~Abhyankar, T.T.~Moh, A reduction theorem for divergent
power series, \textit{J. Reine Angew. Math.}, \textbf{241} (1970), 27--33.

\bibitem{ALM} B.~Al-Shutnawi {\it et al},  On convergence sets of power series with holomorphic coefficients, preprint.

\bibitem{BT} E.~Bedford, B.A.~Taylor, Plurisubharmonic functions with logarithmic singularities, \textit{Ann.\  Inst.\ Fourier (Grenoble)}, \textbf{38} (1988), no.~4, 133--171.

\bibitem{BT2} E.~Bedford, B.A.~Taylor, A new capacity for plurisubharmonic functions. \textit{Acta Math.}, \textbf{149} (1982), no. 1-2, 1--40.

\bibitem{De} J.~Deny, Sur les infinis d'un potentiel, {\it C.\ R.\ Acad.\ Sci.\ Paris S\'er.\ I Math.}, {\bf 224} (1947), 524--525.

\bibitem{GH} P.~Griffiths, J.~Harris, \textit{Principles of algebraic geometry}, John Wiley \& Sons, New York, 1978.

\bibitem{GZ} V.~Guedj, A.~Zeriahi, Intrinsic capacities on compact K\"ahler manifolds, \textit{J.\ Geom.\ Anal.}, \textbf{15} (2005), 607--639.

\bibitem{Ha} F.~Hartogs, Zur Theorie der analytischen Funktionen mehrerer unabh\"angiger Ver\"anderlichen, Math.\ Ann.\ {\bf 62} (1906), 1--88.

\bibitem{HL} F.R.~Harvey, B.~Lawson, Projective hulls and projective Gelfand transform, \textit{Asian J.\ Math.}, \textbf{10} (2006) 607--646.

\bibitem{Kl} M.~Klimek, \textit{Pluripotential theory}, Clarendon Press, New York, 1991.

\bibitem{Le} P.~Lelong, On a problem of M.A.~Zorn, \textit{Proc. Amer.
Math. Soc.}, \textbf{2} (1951), 11--19.

\bibitem{LM} N.~Levenberg, R.E.~Molzon, Convergence sets of a formal
power series, \textit{Math. Z.}, \textbf{197} (1988), 411--420.

\bibitem{LP} N.~Levenberg, E.~Poletsky, Pluripolar hulls, \textit{Michigan Math.\ J.}, \textbf{46} (1999), no.~1, 151--162.

\bibitem{DT} D.~Ma, T.S.~Neelon, On convergence sets of formal power series,  \textit{J. Complex Analysis and its Synergies}, (2015), 1:4, DOI 10.1186/s40627-015-0004-4.

\bibitem{Pe} R.~P\'erez-Marco, A note on holomorphic extensions, preprint, 2000, 	arXiv:math/0009031.

\bibitem{Rt} T.~Ransford, \textit{Potential theory in the complex plane}, Cambridge University Press, New York, 1995.

\bibitem{Sa1} A.~Sadullaev, An estimate for polynomials on analytic sets, \textit{Math.\ USSR Izvestia}, \textbf{20} (1983), 493--502.

\bibitem{Sa2} A.~Sadullaev, Plurisubharmonic measures and capacities on complex manifolds, \textit{Uspekhi Mat.\ Nauk}, \textbf{36} (1981), 53--105; \textit{translation in Russian Math.\ Surveys} \textbf{36} (1981), 61--119.

\bibitem{Sa} A.~Sathaye, Convergence sets of divergent power series, 
\textit{J. Reine Angew. Math.}, \textbf{283} (1976), 86--98.

\bibitem{Sn} N.~Sibony, Sur la fronti\`{e}re de Shilov des domaines de $\CC^n$, Math.\ Ann.\ {\bf 273} (1985), 115--121.

\bibitem{Sc3} J.~Siciak, Extremal phurisubharmonic functions and capacities in $\CC^n$, {\it Sophia Kokyuroku Math.}, {\bf 14} (1982), Sophia University, Tokyo.

\bibitem{Ze} A.~Zeriahi, Ensembles pluripolaires exceptionnels pour la croissance partielle des fonctions holomorphes, \textit{Ann.\ Polon.\ Math.}, \textbf{50} (1989), 81--91.
\end{thebibliography}
\end{document}